\documentclass[12pt]{article}
 
\usepackage{natbib}

\usepackage{amsfonts}

\usepackage{latexsym}

\usepackage{amssymb}

\usepackage{graphicx}

\usepackage{psfrag}

\newcommand{\D}{\displaystyle}

\newcommand{\ignore}[1]{}

\newcommand{\breath}{\medskip}

\newtheorem{thm}{Theorem}[section]

\newcounter{claimcount}[thm]

\newcounter{subclaimcount}[claimcount]

\newtheorem{prop}[thm]{Proposition}

\newtheorem{lemma}[thm]{Lemma}

\newtheorem{cor}[thm]{Corollary}

\newcommand{\dfn}{\sf\em}

\newcommand{\Theorem}[2]{\begin{thm}{\sf #1}  #2 \end{thm}}

\newcommand{\Proposition}[2]{\begin{prop}{\sf #1}  #2 \end{prop}}

\newcommand{\Lemma}[2]{\begin{lemma}{\sf #1}  #2 \end{lemma}}

\newcommand{\Corollary}[2]{\begin{cor}{\sf #1}  #2 \end{cor}} 

\newcommand{\thmfont}[1]{{\sl #1}}    

\newcommand{\example}[1]{
        \refstepcounter{thm}
                     \begin{list}{}
 			{\setlength{\leftmargin}{0em}
 			\setlength{\rightmargin}{0em}}
        \item {\bf Example \thethm.} #1
                   \hfill$\diamondsuit$  \end{list}  
 			}

\newcommand{\bthmlist}{
 \begin{list}{{\bf (\alph{enumii})}}{\usecounter{enumii}}
 			{\setlength{\leftmargin}{0em}
 			\setlength{\itemsep}{0em}
 			\setlength{\parsep}{0em}
 			\setlength{\rightmargin}{0em}}}

\newcommand{\ethmlist}{\end{list}}

\newcommand{\Claim}[1]{\refstepcounter{claimcount}
 
               {\bf Claim \theclaimcount: \ }\thmfont{ #1}}

\newcommand{\Subclaim}[1]{\refstepcounter{subclaimcount}
 
                {\bf Claim \theclaimcount.\thesubclaimcount: \ }\thmfont{ #1}}

\newcommand{\claim}{\Claim}

\newcommand{\subclaim}{\Subclaim}

\newcommand{\bprf}[1][Proof.]{\begin{list}{}
 			{\setlength{\leftmargin}{1em}
 			\setlength{\rightmargin}{0em}}
                         \item {\em \hspace{-1.4em}  #1  }}

\newcommand{\eprf}{\end{list}}

\newcommand{\bthmprf}{\bprf}

\newcommand{\bclaimprf}{\bprf}

\newcommand{\bsubclaimprf}{\bprf}

\newcommand{\ethmprf}{ \hfill$\Box$ 
 \eprf
 
 \breath
 
 }

\newcommand{\eclaimprf}{ \hfill $\Diamond$~{\scriptsize {\tt Claim~\theclaimcount}}\eprf}

\newcommand{\esubclaimprf}{ \hfill $\triangledown$~{\scriptsize 
 {\tt Claim~\theclaimcount.\thesubclaimcount}}\eprf}

\newcommand{\beq}{\begin{eqnarray*}}

\newcommand{\eeq}{\end{eqnarray*}}

\newcommand{\beqn}{ \begin{equation} }

\newcommand{\eeqn}{ \end{equation} }

\newcommand{\If}{\mbox{\ if \ }} 

\newcommand{\And}{\mbox{\ and \ }} 

\newcommand{\Holder}{H\"older}

\newcommand{\dB}{{\mathbb{B}}}

\newcommand{\dF}{{\mathbb{F}}}

\newcommand{\dG}{{\mathbb{G}}}

\newcommand{\dM}{{\mathbb{M}}}

\newcommand{\dN}{{\mathbb{N}}}

\newcommand{\dR}{{\mathbb{R}}}

\newcommand{\dU}{{\mathbb{U}}}

\newcommand{\dV}{{\mathbb{V}}}

\newcommand{\dW}{{\mathbb{W}}}

\newcommand{\dX}{{\mathbb{X}}}

\newcommand{\dZ}{{\mathbb{Z}}}

\newcommand{\barh}{{\overline{h}}}

\newcommand{\bardel }{{\overline{\delta}}}

\newcommand{\barsO}{{\overline{\mathcal{ O}}}}

\newcommand{\bO}{{\mathbf{ O}}}

\newcommand{\ba}{{\mathbf{ a}}}

\newcommand{\bb}{{\mathbf{ b}}}

\newcommand{\bc}{{\mathbf{ c}}}

\newcommand{\bm}{{\mathbf{ m}}}

\newcommand{\bw}{{\mathbf{ w}}}

\newcommand{\bx}{{\mathbf{ x}}}

\newcommand{\by}{{\mathbf{ y}}}

\newcommand{\bz}{{\mathbf{ z}}}

\usepackage{amsbsy}

\newcommand{\sA}{{\mathcal{ A}}}

\newcommand{\sB}{{\mathcal{ B}}}

\newcommand{\sC}{{\mathcal{ C}}}

\newcommand{\sG}{{\mathcal{ G}}}

\newcommand{\sO}{{\mathcal{ O}}}

\newcommand{\sX}{{\mathcal{ X}}}

\newcommand{\sY}{{\mathcal{ Y}}}

\newcommand{\alp }{\alpha}

\newcommand{\bet }{\beta}

\newcommand{\gam }{\gamma}

\newcommand{\del }{\delta}

\newcommand{\eps }{\epsilon}

\newcommand{\lam }{\lambda}

\newcommand{\Gam }{\Gamma}

\newcommand{\fo}{{\mathsf{ o}}}

\newcommand{\fu}{{\mathsf{ u}}}

\newcommand{\fv}{{\mathsf{ v}}}

\newcommand{\fw}{{\mathsf{ w}}}

\newcommand{\fz}{{\mathsf{ z}}}

\newcommand{\tl}{\widetilde}

\newcommand{\tlba}{{\widetilde{\mathbf{ a}}}}

\newcommand{\tlbb}{{\widetilde{\mathbf{ b}}}}

\newcommand{\tlpi }{{\widetilde{\pi}}}

\newcommand{\tlPhi}{{\widetilde{\Phi }}}

\newcommand{\undD}{{\underline{D}}}

\newcommand{\undh}{{\underline{h}}}

\newcommand{\unddel }{{\underline{\delta}}}

\newcommand{\lb}{\left}

\newcommand{\rb}{\right}

\newcommand{\maketall}{\rule[-0.5em]{0em}{1em}}

\newcommand{\implies}{\ensuremath{\Longrightarrow}}

\newcommand{\map}{{\longrightarrow}}

\newcommand{\goto}{{\rightarrow}}

\newcommand{\into}{{\map}}

\newcommand{\seilpmi}{{\Longleftarrow}}

\newcommand{\statement}[1]{\lb(  \maketall 
      \begin{minipage}{40em}
       \begin{tabbing}
         #1 
       \end{tabbing}
      \end{minipage}  \rb)}

\newcommand{\oo}{{\infty}}

\newcommand{\x}{\times}

\newcommand{\union}{\cup}

\newcommand{\Union}{\bigcup}

\newcommand{\intsct}{\cap}

\newcommand{\disj}{\sqcup}

\newcommand{\Disj}{\bigsqcup}

\newcommand{\set}[2]{{\left\{ #1 \; ; \; #2 \right\} }}

\newcommand{\inn}[1]{{\left\langle #1 \right\rangle }}

\newcommand{\choice}[1]{{\lb\{ \begin{array}{rcl}
                                 #1 
                               \end{array}  \rb.  }}

\newcommand{\eeequals}[1]{\raisebox{-0.9ex}{$\overline{\overline{{\scriptscriptstyle{\mathrm{#1}}}}}$}}

\newcommand{\leeeq}[1]{\raisebox{-1ex}{${{\D\leq} \atop {\scriptscriptstyle{\mathrm{#1}}}}$}}

\newcommand{\grt}[1]{\raisebox{-1ex}{${{\D>} \atop {\scriptscriptstyle{\mathrm{#1}}}}$}}

\newcommand{\geeeq}[1]{\raisebox{-1ex}{${{\D\geq} \atop {\scriptscriptstyle{\mathrm{#1}}}}$}}

\newcommand{\subseteeeq}[1]{\raisebox{-1.3ex}{$\stackrel{\D\subseteq}{\scriptscriptstyle{\mathrm{#1}}}$}}

\newcommand{\iiimplies}[1]{\lefteqn{\eeequals{\ #1}}\implies}

\newcommand{\iiiff}[1]{\Leftarrow\!\!\!\!\lefteqn{\eeequals{#1}}\Rightarrow}

\newcommand{\Real}{\dR}

\newcommand{\Natur}{\dN}

\newcommand{\Zahl}{\dZ}

\newcommand{\Zahlmod}[1]{{\Zahl_{/#1}}}

\newcommand{\CC}[1]{{\lb[ #1 \rb]}}

\newcommand{\CO}[1]{{\lb[ #1 \rb)}}

\newcommand{\OO}[1]{{\lb( #1 \rb)}}

\newcommand{\tldV}{{\widetilde{\dV}}}

\newcommand{\AV}{\sA^\dV}

\newcommand{\BW}{\sB^\dW}

\newcommand{\AU}{\sA^\dU}

\newcommand{\AN}{\sA^\Natur}

\newcommand{\BN}{\sB^\Natur}

\newcommand{\AZ}{\sA^\Zahl}

\newcommand{\ZD}{\Zahl^D}

\newcommand{\AZD}{\sA^{\ZD}}

\newcommand{\ZDNE}{\Zahl^D\x\Natur^E}

\newcommand{\AZDNE}{\sA^{\ZDNE}}

\newcommand{\In}{{\!\scriptscriptstyle{\mathrm{in}}}}

\newcommand{\upstream}{\leadsto}

\newcommand{\samestream}{\leftrightsquigarrow}

\newcommand{\connect}{\,\lefteqn{\bullet}\!\!\rightarrow\!}

\newcommand{\boxdim}{\mathrm{boxdim}}

\newcommand{\htop}{h_{\scriptscriptstyle{\mathrm{top}}}}

\newcommand{\bardim}{\overline{\mathrm{dim}}}

\newcommand{\unddim}{\underline{\mathrm{dim}}}

\newcommand{\Remark}[1]{\refstepcounter{thm}\paragraph{Remark \thethm.} #1 \hfill$\diamondsuit$ 
 
 \breath}

\begin{document}
\title{Positive expansiveness versus network dimension in symbolic dynamical systems}
\author{Marcus Pivato \\ Department of Mathematics, Trent University, Canada}
\maketitle

\begin{abstract}
  A `symbolic dynamical system' is a continuous transformation
$\Phi:\sX\into\sX$ of closed perfect subset $\sX\subseteq\AV$, where
$\sA$ is a finite set and $\dV$ is countable. (Examples include
subshifts, odometers, cellular automata, and automaton networks.) \
The function $\Phi$ induces a directed graph structure on $\dV$, whose
geometry reveals information about the dynamical system $(\sX,\Phi)$.
The `dimension' $\dim(\dV)$ is an exponent describing the growth rate
of balls in the digraph as a function of their radius.  We show: if
$\sX$ has positive entropy and $\dim(\dV)>1$, and the system
$(\AV,\sX,\Phi)$ satisfies minimal symmetry and mixing conditions,
then $(\sX,\Phi)$ cannot be positively expansive; this generalizes a
well-known result of Shereshevsky about multidimensional cellular
automata.  We also construct a counterexample to a version of this
result without the symmetry condition.  Finally, we show that network
dimension is invariant under topological conjugacies which are
\Holder-continuous.
\end{abstract}

  Let $\sX$ be {\dfn Cantor space} (the compact, perfect,
zero-dimensional metrizable topological space, which is unique up to
homeomorphism).  A {\dfn Cantor dynamical system} is a continuous
self-map $\Phi:\sX\into\sX$.  In addition to its intrinsic interest,
the class of Cantor systems is important because it has two universal
properties.  First, any topological dynamical system on a compact
metric space is a factor of a Cantor system; see \cite[Corollary 3.9,
p.106]{KurkaBook} or \cite[p.1241]{BalcarSimon}.  Second, the
Jewet-Krieger Theorem says that any ergodic measure-preserving system
can be represented as a uniquely ergodic, minimal Cantor system
\cite[\S4.4, p.188]{Petersen}.

  If $\sA$ is a finite set, and $\dV$ is a countably infinite set,
then the product space $\AV$ is a Cantor space.  Thus, any Cantor
dynamical system can be represented as a self-map $\Phi:\AV\into\AV$,
or more generally, as a self-map $\Phi:\sX\into\sX$, where
$\sX\subset\AV$ is a {\dfn pattern space} (a closed perfect subset of
$\AV$).  We refer to the structure $(\AV,\sX,\Phi)$ as a {\dfn
symbolic dynamical system}.  At an abstract topological level, any
pattern space $\sX$ is homeomorphic to Cantor space, so a symbolic
dynamical system is simply a Cantor dynamical system.  What
distinguishes symbolic dynamical systems is a particular way of
representing $\sX$ as a subset of some Cartesian product $\AV$ (so
that an element of $\sX$ corresponds to some $\dV$-indexed `pattern'
of `symbols' in the alphabet $\sA$).

 The {\dfn network} of $\Phi$ is the digraph structure $(\connect)$ on
$\dV$ defined as follows: for all $\fv,\fw\in\dV$, we have
$\fv\connect\fw$ if and only if the value of $\Phi(\bx)_\fw$ depends
nontrivially on the value of $x_\fv$.  We say that $(\dV,\connect)$
has {\dfn dimension} $\del$ if the cardinality of a ball of radius $r$
grows like $r^\del$ as $r\goto\oo$.  (Note that $\del$ is not
necessarily an integer.) \ For example, if $\Phi:\AZD\into\AZD$ is a
cellular automaton, then its network is just a Cayley digraph on
$\ZD$; the dimension of this network is $D$.  

  This paper explores the relationship between network dimension and
the properties of $(\sX,\Phi)$ as a topological dynamical system.  In
\S\ref{S:prelim}, we formally define the dimension of a network
$(\dV,\connect)$ and the entropy of a pattern space on $\dV$.  In
\S\ref{S:posexp}, we generalize an important result of Shereshevsky
(later reproved by Finelli, Manzini, and Margara) about
multidimensional cellular automata.  We show: if
$\dim(\dV,\connect)>1$ (more generally, if $(\dV,\connect)$ has
`superlinear connectivity'), and $\sX$ has positive entropy and a mild
`mixing' condition, and the system $(\AV,\sX,\Phi)$ has some minimal
symmetries, then $(\sX,\Phi)$ cannot be positively expansive (Theorem
\ref{not.posexpansive}).
In \S\ref{S:propagation}, we consider the {\dfn propagation}
of a symbolic dynamical system, and its relationship with sensitivity
and equicontinuity.   In \S\ref{S:counterexample}, we show that a
`naive' generalization of Shereshevsky's result cannot be true, by
constructing a positively expansive symbolic dynamical system with
network dimension two.  Thus, any result similar to Theorem
\ref{not.posexpansive} must impose at least some additional technical
conditions.

The counterexample in \S\ref{S:counterexample} also shows that network
dimension is not invariant under topological conjugacy; thus, it
cannot be treated as a structural property of an abstract Cantor
dynamical system $(\sX,\Phi)$.  However, in \S\ref{S:lipschitz}, we
propose to augment the system $(\sX,\Phi)$ with a metric which is
Lipschitz for $\Phi$; we show that network dimension {\em is} a
structural property of this `metric' Cantor system, as it is invariant
under all bi\Holder\ conjugacies (Corollary
\ref{network.dimension.invariant.under.holder.conjugacy}).
Sections \ref{S:posexp}-\ref{S:lipschitz} are logically independent, and
can be read in any order.

\section{Preliminaries \label{S:prelim}}

Let $\sA$ be a finite set (called an {\dfn alphabet}) endowed with the
discrete topology.  Let $\dV$ be a countably infinite set of points
(called {\dfn vertices}).  Endow $\AV$ with the Tychonoff product
topology.  A {\dfn pattern space} is a closed perfect subset
$\sX\subseteq\AV$.  A {\dfn symbolic dynamical system} is triple
$(\AV,\sX,\Phi)$, where $\sX\subseteq\AV$ is a pattern space and
$\Phi:\AV\into\AV$ is a continuous function such that
$\Phi(\sX)\subseteq\sX$.  (Sometimes we will simply indicate this as
$(\sX,\Phi)$ when $\sA$ and $\dV$ are clear from context.)

\example{
(a) Let $\dV=\ZDNE$ for some $D,E\geq 0$ and let $\sX\subset\AZDNE$ be
a {\dfn subshift} (i.e. closed, shift-invariant subset).  Then $\sX$
is a pattern space.  Fix $\fz\in\ZDNE$, and let
$\sigma^\fz:\AZDNE\into\AZDNE$ be the associated shift map.  Then
$(\AZDNE,\sX,\sigma^\fz)$ is a symbolic dynamical system.

(b) Let $\dV=\ZDNE$, and let $\Phi:\AZDNE\into\AZDNE$ be a {\dfn
cellular automaton} (CA) ---i.e. a continuous, shift-commuting map.
Then $(\AZDNE,\Phi)$ is a symbolic dynamical system.  More generally, if
$\sX\subset\AZDNE$ is any $\Phi$-invariant subshift, then $(\AZDNE,\sX,\Phi)$
is a symbolic dynamical system.

(c) For all $\fv\in\dV$, let $\sA_\fv\subseteq\sA$. Let
$\sX:=\prod_{\fv\in\dV} \sA_\fv$;  then $\sX$ is a pattern space.
If $\Phi:\sX\into\sX$ is a continuous self-map, then
$(\AV,\sX,\Phi)$ is a symbolic dynamical system,
sometimes called an {\dfn automaton network},
because it can be interpreted as an infinite network of interacting
finite-state automata.

(d) Gromov has initiated a study of `proalgebraic' dynamical systems,
which are (loosely speaking) projective limits of polynomial
self-mappings of algebraic varieties \cite{Gromov}.
If the base field $\dF$ is  finite
(e.g. $\dF=\Zahlmod{p}$), then a `proalgebraic space' can be represented
as a pattern space with alphabet $\dF$; hence a proalgebraic 
system is a symbolic dynamical system.
}

The analysis of subshifts and cellular automata depends heavily on the
highly symmetric structure created by shift-invariance.  Likewise,
Gromov's analysis of proalgebraic systems requires a structure of
`local' symmetries (called {\em holonomies}).  We will also make use
of some minimal symmetry assumptions in \S\ref{S:posexp}.  However, in
general, symbolic dynamical systems do not have any appreciable
symmetries.

For any $\bx\in\sX$ and
$\dU\subset\dV$, we define $\bx_\dU:=[x_\fu]_{\fu\in\dU}\in\sA^\dU$;
we then define $\sX_\dU:=\set{\bx_\dU}{\bx\in\sX}\subseteq\sA^\dU$.
A function $\gam:\AU\into\sB$ is {\dfn proper} if $\gam$ depends nontrivially
on every coordinate in $\dU$.  Formally:  for every $\fv\in\dU$, there
exist $\ba,\ba'\in\AU$ such that: (1) $a_\fu=a'_\fu$ for all $\fu\in\dU\setminus\{\fv\}$; \ (2) $a_\fv\neq a'_\fv$; and (3) $\gam(\ba)\neq \gam(\ba')$.

\Lemma{\label{nhood.corollary}}
{
Let $(\AV,\sX,\Phi)$ be a symbolic dynamical system.
Then for all
$\fv\in\dV$, there is some finite subset $\dU(\fv)\subset\dV$
and a proper function $\phi_\fv:\sA^{\dU(\fv)}\into\sA$ such that,
for any $\bx\in\sX$,  \ $\Phi(\bx)_\fv = \phi_\fv(\bx_{\dU(\fv)})$.
}
\bthmprf
 For each $\fv\in\dV$,  the existence of a local rule
$\phi_\fv$ is proved by exactly the same argument as the Curtis-Hedlund-Lyndon
theorem for cellular automata (see e.g. \cite[Theorem 5.2,
p.190]{KurkaBook}; observe that the construction of the local rule does
not depend on shift-invariance).  The difference is that there
may be different local rules at different vertices.
\ethmprf

The set $\dU(\fv)$ is called the {\dfn input neighbourhood} of $\Phi$
at $\fv$, and denoted $\Phi_\In(\fv)$.  The function
$\phi_\fv:\sA^{\Phi_\In(\fv)}\into\sA$ is called the {\dfn local rule}
of $\Phi$ at $\fw$. 

\example{\label{X:local.rules}
(a) If $\Phi:\AZD\into\AZD$ is a cellular automaton, then
Corollary \ref{nhood.corollary} plus shift-invariance
yields the Curtis-Hedlund-Lyndon theorem.

(b)  Fix $A\in\Natur$ and let $\sA:=\CO{0...A}$.
Let $\bm:=(m_0,m_1,m_2,\ldots)$ be a sequence of natural
numbers in $\CC{1...A}$.
Let $\dV:=\Natur$, and let 
$\sX:=\{\ba\in\AN$;\   $0\leq a_\fv< m_\fv, \ \forall \ \fv\in\dV\}$.
Let $\Phi:\sX\into\sX$ be the $\bm$-ary odometer \cite[\S4.1, p.136]{KurkaBook}.
Then $\Phi_\In(0)=\{0\}$ and $\phi_0:\sA\into\sA$ is
defined by $\phi_0(a_0):=(a_0+1)\bmod{m_0}$.  Meanwhile,
for all $N\geq 1$, we have $\Phi_\In(N)=\CC{0...N}$,
and $\phi_N:\sA^\CC{0...N}\into\sA$ is defined by

$\phi_N(a_0,a_1,\ldots,a_N)\quad:=\quad\choice{(a_N+1) \bmod m_N &&\If
a_n=m_n-1, \ \forall n\in\CO{0...N};\\
a_N &&\mbox{otherwise.}}$
}

\paragraph{Directed graphs.}
Let $\dV$ be a set of `vertices'. A {\dfn directed graph} (or {\dfn
digraph}) structure on $\dV$ is a binary relation $(\connect)$ on
$\dV$ (i.e. a subset $(\connect)\subseteq\dV\x\dV$).  For any
$\fv,\fw\in\dV$, we write ``$\fv\connect\fw$'' if
$(\fv,\fw)\in(\connect)$.  More generally, we say $\fv$ is {\dfn
upstream} of $\fw$ (``$\fv\upstream\fw$'') if either $\fv=\fw$, or
there is a directed path
$\fv=\fv_0\connect\fv_1\connect\cdots\connect\fv_n=\fw$.  The relation
$(\upstream)$ is a partial order (it is reflexive and transitive).  We
write $\fv\samestream\fw$ if $\fv\upstream \fw$ and $\fw\upstream
\fv$.  Thus, $(\samestream)$ is an equivalence relation; the
$(\samestream)$-equivalence classes of $\dV$ are called the {\dfn
biconnected components} of $(\dV,\connect)$.  We say that
$(\dV,\connect)$ is {\dfn biconnected} if $\fv\samestream\fw$ for all
$\fv,\fw\in\dV$. 

  Let $(\sim)$ be the smallest equivalence relation on $\dV$ which
contains $(\connect)$.  Equivalently, for any $\fu,\fw\in\dV$, we have
$\fu\sim\fw$ if either (1) $\fu\upstream\fw$; or (2)
$\fw\upstream\fu$; or (3) (inductively) there exists some $\fv\in\dV$
such that $\fu\sim\fv\sim\fw$.  The $(\sim)$-equivalence classes are
the {\dfn connected components} of $\dV$; if $\fv\sim\fw$ for all
$\fv,\fw\in\dV$, then we say that $(\dV,\connect)$ is {\dfn
connected}.

  If $\Phi:\AV\into\AV$ is any continuous function, then we can define
a digraph relation $(\connect)$ on $\dV$ by $\statement{$\fv\connect
\fw$}\Leftrightarrow \statement{$\fv\in\Phi_\In(\fw)$}$, where
$\Phi_\In(\fw)$ is as defined by Lemma \ref{nhood.corollary} above.
This digraph is called the {\dfn network} of $\Phi$.

\begin{figure}
\centerline{\includegraphics[scale=0.85,clip=true,trim=0 70 0 70]{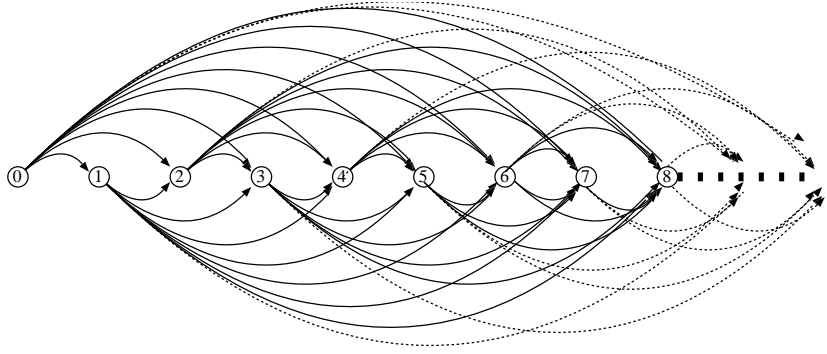}}
\caption{The network of an odometer.\label{fig:odometer.network}}
\end{figure}

\example{\label{X:networks} (a)
Let $\Phi:\AZDNE\into\AZDNE$ be a cellular automaton; then the network
of $\Phi$ is a Cayley digraph on $\ZDNE$. 

(b) Figure \ref{fig:odometer.network} depicts the network of the
odometer $\Phi:\AN\into\AN$ from Example \ref{X:local.rules}(b).}

\paragraph{Connectivity dimension.}
Let $(\dV,\connect)$ be an infinite digraph (e.g. the network of
a continuous function $\Phi:\AV\into\AV$).
For any subset $\dU\subset\dV$, define
$\dB(\dU,1):=\dU\union\set{\fv\in\dV}{ \exists \ \fu\in\dU: \
\fv\connect\fu}$.   Then inductively define
$\dB(\dU,n+1):=\dB[\dB(\dU,n),1]$ for all $n\in\Natur$.  Thus,
$\dB(\fw,1):=\{\fw\}\union \Phi_\In(\fw)$, and
$\dB(\fw,r)$ is the set of all $\fv\in\dV$ such that there exists some
path $\fv= \fv_1\connect\fv_2\connect\cdots\connect\fv_s=\fw$
with $s\leq r$.   For any $\fv\in\dV$, we define
\beqn
\label{connectivity.dimension}
\unddim_\fv(\dV,\connect)\ \ :=\ \  \liminf_{r\goto\oo} \ \frac{\log\lb|\dB(\fv,r)\rb|}{\log(r)}
\ \ \And\ \ 
\bardim_\fv(\dV,\connect)\ \ :=\ \  \limsup_{r\goto\oo} \ \frac{\log\lb|\dB(\fv,r)\rb|}{\log(r)}.
\eeqn
If $\unddim_\fv(\dV,\connect)=\bardim_\fv(\dV,\connect)$,
then we refer to their common value as ``$\dim_\fv(\dV,\connect)$'',
the {\dfn connectivity dimension} of $(\dV,\connect)$ at $\fv$,
and we say that $(\dV,\connect)$ is {\dfn dimensionally regular} at
$\fv$.

\example{\label{X:connectivity.dimension}
For all $r\in\Natur$, let $\bet_\fv(r):=|\dB(\fv,r)|$.

(a) Let $\del\in\CO{0,\oo}$, and suppose 
\[
0 \quad<\quad \liminf_{r\goto\oo} \ \frac{\bet_\fv(r)}{r^\del}
\quad\leq \quad \limsup_{r\goto\oo} \ \frac{\bet_\fv(r)}{r^\del}
\quad<\quad \oo.
\]
(For example, suppose $\bet_\fv(r) = C\, r^\del + p(r)$, where $C$ is a constant
and $p$ is a polynomial of degree less than $\del$.)\  Then
$\dim_\fv(\dV,\connect)=\del$.

(b) Likewise, if $C, \del,\lam > 0$, and
$\bet_\fv(r) = C\, r^\del\cdot \log(r)^\lam$, 
 then $\dim_\fv(\dV,\connect)=\del$.

(c) Let $c>0$.  If $\bet_\fv(r) = c^r$, then $\dim_\fv(\dV,\connect)=\oo$.
}

Let
$\bardim(\dV,\connect):=\sup\set{\bardim_\fv(\dV,\connect)}{\fv\in\dV}$
and
$\unddim(\dV,\connect):=\inf\set{\unddim_\fv(\dV,\connect)}{\fv\in\dV}$.
If $\unddim(\dV,\connect)=\bardim(\dV,\connect)$, then we refer to
their common value as ``$\dim(\dV,\connect)$'', the (global) {\dfn
connectivity dimension} of $(\dV,\connect)$, 
and we say that $(\dV,\connect)$ is {\dfn dimensionally homogeneous}.
(This implies that  $(\dV,\connect)$ is everywhere dimensionally regular.)

\example{(a) If $\Zahl^D$ has the obvious Cayley digraph structure,
then $\dim(\Zahl^D)=D$.

(b) Let $(\Natur,\connect)$ be the digraph in Figure \ref{fig:odometer.network}
(the odometer). 
Then $\dim(\Natur,\connect)=0$, because
for all $n\in\Natur$, and all $r\geq 1$, we have $|\dB(n,r)|=n+1$
(because $\dB(n,r)=\CC{0...n}$).

(c) If $(\dV,\connect)$ is the Cayley digraph of a group $\sG$, then
$\dim(\dV,\connect)$ is the `growth dimension' of $\sG$.
More generally,
if $(\dV,\connect)$ is any graph whose automorphism group acts
transitively, then $(\dV,\connect)$ is `almost' a Cayley digraph;
for a survey of the well-developed dimension theory for such graphs, see
\cite[\S5]{MoharWoessSpectra} or
\cite{ImrichSeifterGraphDimSurvey}.}

  The dimension of a Cayley digraph is always an integer.  However, there
exist `self-similar' graphs with fractional
connectivity dimensions \cite{McCannaGraphDim}.  Connectivity
dimension is closely related to properties of diffusion processes and
electrical conductance on graphs
\cite{Telcs1,Telcs2,Telcs3,Telcs4}, the existence of periodic points
in `majority vote' networks \cite{Moran,GinosarHolzman}, and also
arises in certain models of quantum gravity
\cite{NowotnyRequardt1,NowotnyRequardt2}.

   Not all digraphs are dimensionally regular.  For example, consider
a digraph which consists of increasingly large `clumps' which are
spaced at increasingly long intervals along an infinite line-graph; by
making the clumps and the intervals between them grow fast enough, one
can force $\unddim_\fv(\dV,\connect)=1$ while
$\bardim_\fv(\dV,\connect)>1$ for some $\fv\in\dV$.  (However,
examples like this are highly contrived; probably, most `natural'
examples are dimensionally regular.)

Furthermore, not all connected, dimensionally regular digraphs are
dimensionally homogeneous. For example, let $\dV_1\cong\Zahl$ be a
biconnected Cayley digraph of $\Zahl$, and let $\dV_2\cong\Zahl^2$ be
a biconnected Cayley digraph of $\Zahl^2$.  Let $\dV=\dV_1\disj\dV_2$,
with connections $n\connect (n,0)$ for all $n\in \Zahl\cong\dV_1$.
Then $\dV_1$ and $\dV_2$ are biconnected components of $\dV$, with
$\dV_1$ upstream from $\dV_2$.  Clearly, $\dim_\fv(\dV_k)=k$ for all
$\fv\in\dV_k$ and $k=1,2$.

\Lemma{\label{constant.nhood.growth.rate}}
{
 Let $(\dV,\connect)$ be a digraph.  If $\fv\upstream\fw$,
then $\unddim_\fv(\dV,\connect)\leq \unddim_\fw(\dV,\connect)$ and
$\bardim_\fv(\dV,\connect)\leq \bardim_\fw(\dV,\connect)$.
}
\bthmprf
 If $\fv\upstream\fw$, then there exists $R>0$ such that $\fv\in\dB(\fw,R)$.
Thus, for all $r\in\Natur$, we have
$\dB(\fv,r)\subseteq\dB(\fw,R+r)$, hence
$|\dB(\fv,r)|\leq |\dB(\fw,R+r)|$.
  Thus
\beq
\unddim_\fv(\dV,\connect)
&:=&
\liminf_{r\goto\oo} \ \frac{\log\lb|\dB(\fv,r)\rb|}{\log(r)}
\quad\leq\quad
\liminf_{r\goto\oo} \ \frac{\log\lb|\dB(\fw,R+r)\rb|}{\log(R+r)}\cdot \frac{\log(R+r)}{\log(r)}
\\&=&
\lb(\liminf_{r\goto\oo} \ \frac{\log\lb|\dB(\fw,R+r)\rb|}{\log(R+r)}\rb)
\cdot \lb(\lim_{r\goto\oo} \frac{\log(R+r)}{\log(r)}\rb)
\quad=\quad \unddim_\fw(\dV,\connect)\cdot 1.
\eeq
Hence  $\unddim_\fv(\dV,\connect)\leq \unddim_\fw(\dV,\connect)$.  Likewise,
$\bardim_\fv(\dV,\connect)\leq \bardim_\fw(\dV,\connect)$.
\ethmprf

 If $\dW\subset\dV$ is a biconnected component of
$(\dV,\connect)$, then Lemma \ref{constant.nhood.growth.rate} says
that every vertex in $\dW$ has the same connectivity dimension.  In particular, if $(\dV,\connect)$
is biconnected and dimensionally regular, then
it is dimensionally homogeneous.

\paragraph{Entropy.}
Let $(\dV,\connect)$ be a digraph, and let $\sX\subset\AV$ be a
pattern space.   For any
 $\fv\in\dV$, we define the {\dfn lower} and {\dfn upper topological
 entropies} of $\sX$ around $\fv$ by: \beqn
\label{entropy.defn}
  \undh_\fv(\sX)\quad:=\quad \liminf_{r\goto\oo} \frac{\log_2|\sX_{\dB(\fv,r)}|}
{|\dB(\fv,r)|}
\quad\And\quad
  \barh_\fv(\sX)\quad:=\quad \limsup_{r\goto\oo} \frac{\log_2|\sX_{\dB(\fv,r)}|}
{|\dB(\fv,r)|}.
\eeqn
Clearly, $0\leq \undh_\fv(\sX)\leq\barh_\fv(\sX)\leq \log_2|\sA|$.
Let $\undh(\sX):=\D \inf_{\fv\in\dV}\ h_\fv(\sX)$ and
$\barh(\sX):=\D \sup_{\fv\in\dV}\ h_\fv(\sX)$.

\section{Positive expansion versus network connectivity
\label{S:posexp}}

  An abstract Cantor dynamical system  $(\sX,\Phi)$ is {\dfn posexpansive}
if it is topologically conjugate to a one-sided shift.
In particular, let $(\AV,\sX,\Phi)$ be a symbolic dynamical system.
Fix a finite subset $\dW\subset\dV$, let $\sB:=\sA^\dW$
and define the function $\Phi^\Natur_\dW:\sX\into\BN$ by
\[
\Phi^\Natur_\dW(\bx) \quad:=\quad[\bx_\dW, \ \Phi(\bx)_\dW, \ \Phi^2(\bx)_\dW,
 \ \Phi^3(\bx)_\dW, \ldots],
\qquad\mbox{for all $\bx\in\sX$.}
\]
Then $(\sX,\Phi)$ is posexpansive if and only if there is some finite subset
$\dW\subset\dV$ (called a {\dfn posexpansive window}) such that the
function $\Phi^\Natur_\dW$ is an injection.  If
$\sY:=\Phi^\Natur_\dW[\sX]\subseteq\BN$ and $\sigma:\BN\into\BN$ is
the shift map, then  $\sigma(\sY)\subseteq\sY$, and 
$\Phi^\Natur_\dW$ is a
topological conjugacy from $(\sX,\Phi)$ to the system $(\sY,\sigma)$.

For any $\dW\subset\dV$ and  $T\in\Natur$, define 
$\Phi^{\CC{0...T}}_\dW:\sX\into\sB^{\CC{0...T}}$ by
\[
\Phi^{\CC{0...T}}_\dW(\bx) \quad:=\quad[\bx_\dW, \ \Phi(\bx)_\dW, \ \Phi^2(\bx)_\dW, \ldots, \Phi^T(\bx)_\dW],
\qquad\mbox{for all $\bx\in\sX$.}
\]
Let $\dW^T:=\set{\fv\in\dV}{\forall \ \bx,\bx'\in\sX, \ 
\statement{$\Phi^{\CC{0...T}}_\dW(\bx)=\Phi^{\CC{0...T}}_\dW(\bx')$}\implies
\statement{$x_\fv=x'_\fv$}}$.  Thus,
\beqn
\label{expansive.panorama1}
\mbox{for all $\bx,\bx'\in\sX$,}\qquad
\statement{$\Phi^{\CC{0...T}}_\dW(\bx)=\Phi^{\CC{0...T}}_\dW(\bx')$}\quad\implies\quad
\statement{$\bx_{\dW^T}=\bx'_{\dW^T}$}.
\eeqn
Then we have
$\dW=\dW^0\subseteq\dW^1\subseteq\dW^2\subseteq\dW^3\subseteq\cdots$;
the sequence $\{\dW_0^t\}_{t=0}^\oo$ is called 
the {\dfn $(\sX,\Phi)$-panorama} of $\dW$.
Clearly, $\dW$ is a posexpansive window for $\Phi$ if and only if
\beqn
\label{expansive.panorama2}
\Union_{t=0}^\oo \dW^t \quad=\quad \dV. 
\eeqn
  Shereshevsky has shown that multidimensional cellular automata
can never be posexpansive.  To be precise, he showed: if $(\dG,\cdot)$ is
any group with growth dimension $D\geq2$ (e.g. $\dG=\Zahl^D$), and
$\sX\subset\sA^\dG$ is a subshift with positive topological entropy,
and $\Phi:\sA^\dG\into\sA^\dG$ is an $\sX$-preserving cellular automaton,
then the system $(\sX,\Phi)$ is not posexpansive; see
\cite[Corollary 2]{Shereshevsky93} or
 \cite[Theorem 1.1]{Shereshevsky96}.
The special case when $\dG=\ZD$ and $\sX=\AZD$
was later reproved in 
\cite[Theorem  4.4]{FinelliManziniMargara}.
In this section, we will generalize this result to any symbolic
dynamical system satisfying some mild symmetry and mixing conditions.

\breath

Let $(\dV,\connect)$ be a digraph.  A injection $\tau:\dV\into\dV$ is
a {\dfn subisometry} if, for all $\fv,\fw\in\dV$,
$(\fv\connect\fw)\iff(\tau(\fv)\connect\tau(\fw))$.  Thus, for all
$\fv\in\dV$ and $r>0$, we have $\tau[\dB(\fv,r)]\subseteq \dB[\tau(\fv),r]$ (with equality
if $\tau:\dV\into\dV$ is surjective).  The map
$\tau$ induces a surjection $\tau_*:\AV\into\AV$ defined by $\tau_*(\ba):=\ba'$
where $a'_{\fv} := a_{\tau(\fv)}$ for all $\fv\in\dV$.  Let
$\sX\subseteq\AV$ be a pattern space; if $\tau$ is a subisometry and
$\tau_*(\sX)=\sX$, then we say $\tau$ is a {\dfn subsymmetry} of $\sX$.

\example{Let  $\dV=\ZDNE$ for some $D,E\geq 0$
(or some other finitely generated monoid),
with the Cayley digraph structure induced by some finite generating set.
Fix $\fw\in\dV$, and define the {\dfn shift map}
$\tau^\fw:\dV\into\dV$ by $\tau^\fw(\fv):=\fv+\fw$;
then $\tau$ is a subisometry of the Cayley digraph.
If $\sX\subseteq\AZDNE$ is a
subshift, then $\tau$ is a subsymmetry of $\sX$.  

{\em Note.} Subsymmetries of $\sX$ are not necessarily 
injective.  For example, the unilateral shift on $\AN$ is a subsymmetry,
but it is many-to-one.
}

Let $(\AV,\sX,\Phi)$ be a symbolic dynamical system.
A {\dfn subsymmetry} of 
$(\AV,\sX,\Phi)$ is a subisometry $\tau:\dV\into\dV$ such that
$\tau_*[\sX]=\sX$ and $\tau_*\circ\Phi=\Phi\circ\tau_*$.

\example{\label{X:subsymmetry}
(a) Let $\dV=\ZDNE$ (or any other finitely generated group), let
$\sX\subseteq\AZDNE$ be a subshift, and let $\Phi:\AZDNE\into\AZDNE$ be
an $\sX$-preserving cellular automaton. Then any $(\ZDNE)$-shift
is a subsymmetry of the system $(\AZDNE,\sX,\Phi)$.  

\begin{figure}
\centerline{\includegraphics[scale=0.25,angle=-90]{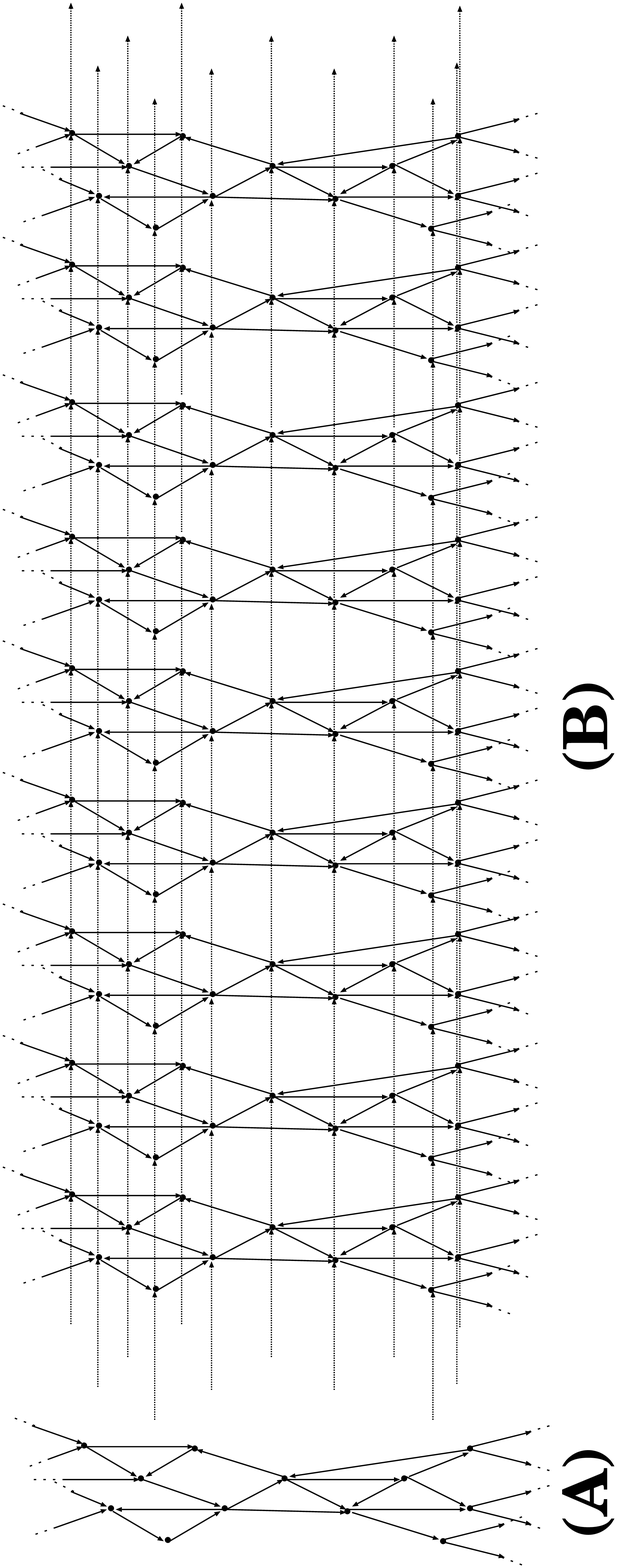}}
\caption{Example \ref{X:subsymmetry}(b).\label{fig:subsymmetry}}
\end{figure}

(b) Let $\psi:\sA\x\sA\into\sA$ be a binary operator (e.g. a group
operator).  Let $(\AV,\Phi)$ be an arbitrary symbolic dynamical
system (perhaps with no symmetries), such as the one in Figure \ref{fig:subsymmetry}(A).  Define $\tldV:=\dV\x\Zahl$, and
identify $\sA^{\tldV}$ with $(\AV)^\Zahl$ in
the obvious way; a generic element of $\sA^{\tldV}$ could be indicated as
$\tlba:=[\ba^n]_{n\in\Zahl}$, where $\ba^n\in\AV$ for all $n\in\Zahl$.
Let $\sigma:(\AV)^\Zahl\into(\AV)^\Zahl$ be the shift map; then $\sigma$
is a subsymmetry of $\sA^\tldV$.  Define $\tlPhi:\sA^\tldV\into\sA^\tldV$ by $\tlPhi[\tlba]
= \tlbb$, where $b^n_\fv = \psi[\Phi(\ba^n)_\fv, \ a^{n+1}_\fv]$ for all
$n\in\Zahl$ and $\fv\in\dV$; this yields the connection network
in Figure \ref{fig:subsymmetry}(B).  Then $\sigma$ is a subsymmetry of 
$(\sA^{\tldV},\tlPhi)$.  }

We say that the pattern space $\sX$ has 
{\dfn weak independence} if there
is some constant $\eps>0$ such that,
for any disjoint balls $\dB_1,\ldots,\dB_N\subset\dV$,
\beqn
\label{long.range.correlations}
  \log_2\lb|\sX_{\dB_1\disj\cdots\disj\dB_N}\rb|
\quad\geq\quad \eps \sum_{n=1}^N \log_2\lb|\sX_{\dB_n}\rb|.
\eeqn
This can be seen as a kind of `topological mixing' condition ---it means
that the information contained in balls $\dB_1,\ldots,\dB_{N-1}$ has limited
power to predict the contents of ball $\dB_N$.

\example{For all $\fv\in\sA$, let $\sA_\fv\subset\sA$ be a subset of
cardinality at least 2.  Let $\sX:=\prod_{\fv\in\dV} \sA_\fv\subset\AV$;
then $\undh(\sX)\geq 1$, and $\sX$ has weak independence.

In particular, the space $\sX=\AV$ itself satisfies weak independence.}

\newcommand{\speed}{\mathrm{speed}}

\ignore{:=\min\set{r+s}{r,s\in\Natur \And
\dB(\fv,r)\intsct\dB(\fw,s)\neq\emptyset}}

For any $\fv\sim\fw\in\dV$, let $d(\fv,\fw)$ be the length of
the shortest {\em undirected} path from $\fv$ to $\fw$;
then $d$ is a metric on each connected component of $\dV$.
(If $\fv\not\sim\fw$, let $d(\fv,\fw):=\oo$).
For any $\fv\in\dV$, let
$\D \speed(\fv,\tau) \ :=\ \D \lim_{n\goto\oo} \frac{d[\fv,\tau^n(\fv)]}{n}$.

\Lemma{\label{speed.lemma}}{
\bthmlist
\item For any $\fv,\fw\in\dV$, we have $d\lb[\tau(\fv),\tau(\fw)\rb] \leq d(\fv,\fw)$ 
{\rm(with equality if $\tau:\dV\into\dV$ is surjective).}

\item For any $\fv\in\dV$, we have
$\speed(\fv,\tau) \D = \inf_{n\in\Natur} \ \frac{d[\fv,\tau^n(\fv)]}{n}$.

\item For all $\fv,\fw\in\dV$, if $\fv\sim\fw$, then
 $\speed(\fv,\tau)=\speed(\fw,\tau)$.
\ethmlist
}
\bthmprf (a) Let $(\fv_0,\fv_1,\fv_2,\ldots,\fv_N)$ be a minimal undirected path
from $\fv$ to $\fw$ (i.e. $\fv_0=\fv$, $\fv_N=\fw$, and either $\fv_{n-1}\connect \fv_{n}$ or
$\fv_{n}\connect \fv_{n-1}$ for all $n\in\CC{1...N}$).  Then
$(\tau(\fv_0),\tau(\fv_1),\ldots,\tau(\fv_N))$ is an undirected path of length $N$
from $\tau(\fv)$ to $\tau(\fw)$.  (However, there may exist shorter paths from
 $\tau(\fv)$ to $\tau(\fw)$ which do not arise as $\tau$-images of paths from $\fv$ to $\fw$).

(b)  is because the sequence
$\{d[\fv,\tau^n(\fv)]\}_{n=1}^\oo$ is subadditive:
\[
d[\fv,\tau^{n+m}(\fv)]\quad\leeeq{(\triangle)}\quad 
d[\fv,\tau^{n}(\fv)]+d[\tau^n(\fv),\tau^{n+m}(\fv)]
\quad\leeeq{(@)}\quad d[\fv,\tau^{n}(\fv)]+d[\fv,\tau^{m}(\fv)].
\]
Here $(\triangle)$ is the triangle inequality, and $(@)$ is by part (a).

To see (c), let $r:=d(\fv,\fw)$ (finite because $\fv\sim\fw$).
Then for any $n\in\Natur$, 
\beq
d[\fv,\tau^n(\fv)]
&\leeeq{(\triangle)}&
d[\fv,\fw] + d[\fw,\tau^n(\fw)]+ d[\tau^n(\fw),\tau^n(\fv)]
\\&\leeeq{(@)}&
d[\fv,\fw] + d[\fw,\tau^n(\fw)]+ d[\fw,\fv]
\quad\eeequals{(\dagger)}\quad
d[\fw,\tau^n(\fw)]+ 2 r.\\
\mbox{Thus,} \  \
\speed(\fv,\tau) &:=& \lim_{n\goto\oo} \frac{d[\fv,\tau^n(\fv)]}{n}
\quad\leq\quad
\lim_{n\goto\oo} \frac{d[\fw,\tau^n(\fw)]+2r}{n}
\quad=\quad \speed(\fw,\tau).
\eeq
Here, $(@)$ is by part (a), 
and $(\dagger)$ is because $d[\fw,\fv]=d[\fv,\fw]$ (because
the definition of $d$ is symmetric) and $d[\fv,\fw]=r$.
A symmetric argument yields $\speed(\fw,\tau)\leq\speed(\fv,\tau)$.
Thus, $\speed(\fv,\tau)=\speed(\fw,\tau)$.
\ethmprf

Lemma \ref{speed.lemma}(b) says the limit defining $\speed(\fv,\tau)$
exists for all $\fv\in\dV$.  We say that $\tau$ is a {\dfn moving
subsymmetry} if $\speed(\fv,\tau)>0$ for all $\fv\in\dV$.  [Lemma
\ref{speed.lemma}(c) implies that it suffices to require
$\speed(\fv,\tau)>0$ for at least one $\fv$ in each connected
component of $(\dV,\connect)$.]

\example{\label{X:speed} (a) \ Let $\dV=\ZD$ with the Cayley digraph
structure induced by the standard generating set $\{(\pm
1,0,0,\ldots,0)$, \ $(0,\pm 1,0,\ldots,0)$, $\ldots, (0,\ldots,0,\pm
1)\}$.  If $\fz=(z_1,\ldots,z_D)\in\ZD$, and $\fo=(0,\ldots,0)$, then
$d(\fo,\fz)=|z_1|+\cdots+|z_D|$.  Let $\sigma^\fz:\ZD\into\ZD$ be the
shift map.  Then $\speed(\sigma^\fz,\fv)=d(\fo,\fz)$ for all
$\fv\in\ZD$.  Thus, all nontrivial shifts are moving symmetries of
$\AV$.

\begin{figure}
\centerline{\includegraphics[scale=0.25,angle=-90,clip=true,trim=0 0 0 4400]{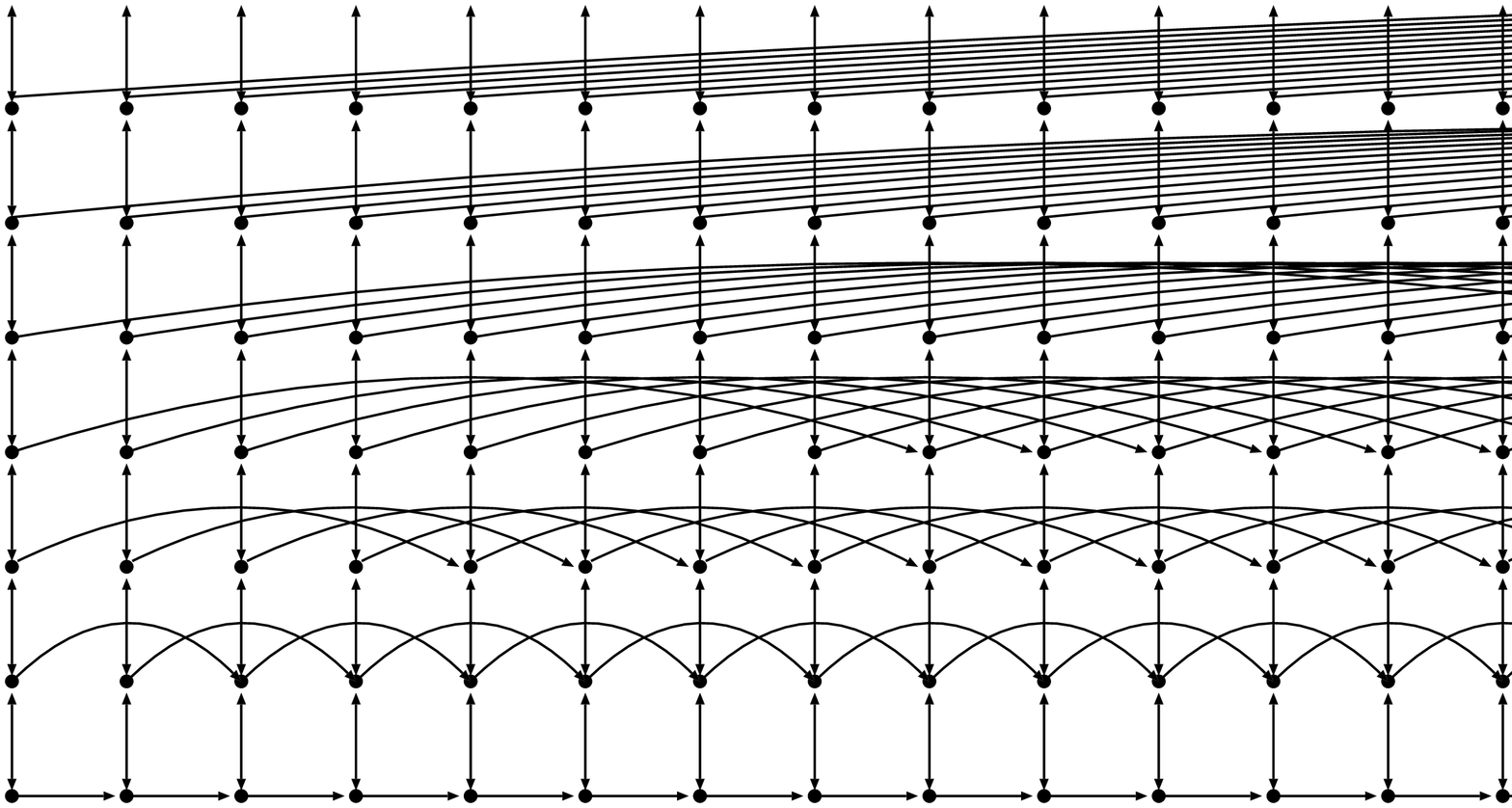}}
\caption{Example \ref{X:speed}(b).\label{fig:speed}}
\end{figure}

(b) \ Let $\dV=\Zahl\x\Natur$, with the digraph structure shown in
Figure \ref{fig:speed}.  Here, for any $(z,n)\in\dV$,
 we have $(z,n)\connect(z,n')$ whenever
$n'=n\pm 1$, and we also have $(z,n)\connect(z',n)$
whenever $z'=z+2^n$.  Thus: $\cdots\connect(-1,0)\connect(0,0)\connect(1,0)\connect(2,0)\connect(3,0)\connect\cdots$, and $\cdots\connect(-2,1)\connect
(0,1)\connect(2,1)\connect(4,1)\connect(6,1)\connect\cdots$,
and $\cdots\connect(-4,2)\connect(0,2)\connect(4,2)\connect(8,2)\connect(12,2)\connect\cdots$, etc.
Define subisometry $\tau:\dV\into\dV$ by $\tau(z,n)=(z+1,n)$.
Then $\speed(\tau,\fv)=0$, for all $\fv\in\dV$, because for any
$k\in\Natur$, there is a
path from $\fv$ to $\tau^{(2^k)}(\fv)$ of length at most $2k+1$.
Thus, $\tau$ is {\em not} a moving subsymmetry.
}

In a digraph $(\dV,\connect)$,
a vertex $\fv\in\dV$  has {\dfn superlinear connectivity}
if $\D\liminf_{r\goto\oo} \ \frac{\lb|\dB(\fv,r)\rb|}{r} \ = \ \oo$.

\example{(a) If 
$\unddim_\fv(\dV,\connect)>1$, then $\fv$
has superlinear connectivity.
(For example, if  $\dV$ is a Cayley digraph
of a group with growth dimension $D\geq 2$, then every
vertex has superlinear connectivity.)

(b) If $\fv\upstream\fw$ and $\fv$ has superlinear connectivity,
then $\fw$ has superlinear connectivity.
}

The main result of this section is the following:

\Theorem{\label{not.posexpansive}}
{
Let $(\AV,\sX,\Phi)$ be a symbolic dynamical system with a moving subsymmetry.
If
$\sX$ has weak independence, and
there exists some $\fv\in\dV$ with superlinear connectivity
such that $\barh_\fv(\sX)>0$, then the system $(\sX,\Phi)$ is not posexpansive. 
}

Before proving  Theorem \ref{not.posexpansive}, we give two
concrete corollaries.

\Corollary{\label{not.posexpansive.corollary1}} 
{ 
Let $(\AV,\sX,\Phi)$ be a symbolic dynamical system with a moving subsymmetry, such that
$\sX$ has weak independence. Suppose that either
\bthmlist
  \item $\barh(\sX)>0$ and $\unddim(\dV,\connect)>1$; or
  \item $(\dV,\connect)$ is dimensionally regular, 
 $\undh(\sX)>0$, and $\bardim(\dV,\connect)>1$.
\ethmlist
Then the system $(\sX,\Phi)$ is not posexpansive. 
}
\bthmprf
Recall: if $\unddim_\fv(\dV,\connect)>1$, then 
$\fv$ has superlinear connectivity.

(a) There exists $\fv\in\dV$ with $\barh_\fv(\sX)>0$,
because $\barh(\sX)>0$.  But $\fv$ also
has superlinear connectivity, because
$\unddim_\fv(\dV,\connect)\geq \unddim(\dV,\connect)>1$.  
Now apply  Theorem \ref{not.posexpansive}.

(b) There exists $\fv\in\dV$ with $\bardim_\fv(\dV,\connect)>1$,
because $\bardim(\dV,\connect)>1$.  Thus,  ${\unddim_\fv(\dV,\connect)}>1$ also,
because $(\dV,\connect)$ is dimensionally regular.
Thus, $\fv$
has superlinear connectivity.
Also,  $\barh_\fv(\sX)>\undh_\fv(\sX)\geq \undh(\sX)>0$.
Now apply Theorem \ref{not.posexpansive}.
\ethmprf

\Corollary{\label{not.posexpansive.corollary2}}
{
Let $\Phi:\AV\into\AV$ be a continuous self-map with a moving
subsymmetry.  If  $\unddim(\dV,\connect)>1$, then the system $(\AV,\Phi)$ is not posexpansive. 
}
\bthmprf
If $\sX=\AV$, then clearly $\sX$ has weak independence, and
$\barh(\sX)=\log_2|\sA|>0$.
Now apply Corollary \ref{not.posexpansive.corollary1}(a).
\ethmprf

The proof of Theorem \ref{not.posexpansive} consists of two lemmas
concerning the `entropy' of a pattern space relative to a subsymmetry.
Let $\sX\subseteq\AV$ be a pattern space and let $\tau:\dV\into\dV$
be a subsymmetry of $\sX$.   For any finite $\dF\subset\dV$, we define
\beqn
\label{tau.entropy.defn1}
  \barh(\sX,\tau,\dF)
\quad:=\quad
\limsup_{N\goto\oo}\ \frac{1}{N} \log_2\lb|\sX_{\dF(N)}\rb|,
\quad
\mbox{where}\quad
\dF(N)\ := \ \Union_{n=0}^N \tau^n(\dF)  \ \subseteq \ \dV.\qquad
\eeqn
We then define the {\dfn upper $\tau$-entropy} of $\sX$ by
\beqn
\label{tau.entropy.defn2}
\barh(\sX,\tau)\quad:=\quad \sup_{{\dF\subset\dV} \atop{\mathrm{finite}}} \ \barh(\sX,\tau,\dF).
\eeqn

\ignore{For any $\dF\subset\dV$, note that
$\barh(\sX,\tau,\dF)\leq N_\dF\cdot \log_2|\sA|$, where $N_\dF$ is the number
of distinct infinite $\tau$-orbits in $\dV$ which $\dF$ intersects.
Thus, $\barh(\sX,\tau)\leq N\cdot \log_2|\sA|$, where $N$ is the
number of infinite $\tau$-orbits in $\dV$.  (In particular,
$\barh(\sX,\tau)$ can be infinite only if $\dV$ contains infinitely many
distinct infinite $\tau$-orbits).}

\Lemma{\label{infinite.tau.entropy}}
{
Let $(\dV,\connect)$ be a digraph,
and let $\sX\subseteq\AV$  be a pattern space with weak independence.
Suppose there exists $\fv\in\dV$ with superlinear connectivity and
$\barh_\fv(\sX)>0$.
If $\tau:\dV\into\dV$ is any moving subsymmetry of $\sX$,
then $\barh(\sX,\tau)=\oo$.
}
\bthmprf Let $\eps>0$ be as in equation (\ref{long.range.correlations}).
Let $S:=\speed(\tau,\fv)>0$.

\Claim{For any $r>0$, we have $\D \barh[\sX,\tau,\dB(\fv,r)]
\ \geq \ 
\frac{S \eps}{4 r} \ \log_2\lb|\sX_{\dB(\fv,r)}\rb|$.
\label{infinite.tau.entropy.C1}}
\bclaimprf
  Let $m:=\lceil 2 r/S\rceil$, then
the points $\{\fv, \ \tau^m(\fv), \  \tau^{2m}(\fv),  \  \tau^{3m}(\fv),\ldots\}$ are all
at least $2r$-separated, by Lemma \ref{speed.lemma}(b).
Thus, the balls $\{\dB(\fv,r)$, \ 
$\dB(\tau^m(\fv),r)$, \  $\dB(\tau^{2m}(\fv),r), \ldots\}$
are all disjoint.  Let $\dF:=\dB(\fv,r)$; 
then for any $n\in\Natur$, we have 
$\tau^{n m}[\dF]\subseteq \dB(\tau^{n m}(\fv),r)$
(because $\tau$ is a subisometry of $\dV$).
Thus, the sets
$\{\dF$, \  $\tau^m(\dF)$, \ $\tau^{2m}(\dF),\ldots\}$ are disjoint.

For any $N\in\Natur$,
let $\dF(Nm)\ := \ \D \Union_{k=0}^{Nm} \tau^k(\dF)$;
 then $\dF(Nm)\ \supseteq \ \D \Disj_{n=0}^N \tau^{nm}(\dF)$.
Thus
\beqn
\label{infinite.tau.entropy.e1}
\log_2\lb|\sX_{\dF(Nm)}\rb| \quad\geeeq{(*)}\quad 
\eps \cdot \sum_{n=1}^N \log_2\lb|\sX_{\tau^{nm}(\dF)}\rb|
\quad\eeequals{(\dagger)}\quad 
\eps N \log_2\lb|\sX_{\dF}\rb|,
\eeqn
where $(*)$ is by equation (\ref{long.range.correlations}),
and $(\dagger)$ is because $\tau$ is a subsymmetry of $\sX$
(so $\lb|\sX_{\dF}\rb| = \lb|\sX_{\tau^{k}(\dF)}\rb|$ for
all $k\in\Zahl$). 
 Combining
equations (\ref{tau.entropy.defn1})
and (\ref{infinite.tau.entropy.e1}), we get
\beq
\barh(\sX,\tau,\dF)
&:\eeequals{(\ref{tau.entropy.defn1})}&
\limsup_{N\goto\oo} \frac{1}{N} \log_2\lb|\sX_{\dF(N)}\rb|
\quad\geq\quad
\limsup_{N\goto\oo} \frac{1}{Nm} \log_2\lb|\sX_{\dF(Nm)}\rb|
\\ &\geeeq{(\ref{infinite.tau.entropy.e1})}&
\lim_{N\goto\oo} \frac{\eps \, N}{N\, m} \log_2\lb|\sX_{\dF}\rb|
\quad=\quad
 \frac{\eps}{m} \cdot \log_2\lb|\sX_{\dF}\rb|
\quad\geeeq{(*)}\quad
 \frac{S\eps}{4 r} \cdot \log_2\lb|\sX_{\dF}\rb|,
\eeq
as desired. 
Here, $(*)$ is because $m:=\lceil 2 r/S\rceil \leq 4 r/S$.
\eclaimprf

It follows from defining equation (\ref{tau.entropy.defn2}) that
\beq
\barh[\sX,\tau]
&\geeeq{(\ref{tau.entropy.defn2})}&
\sup_{r\in\Natur}\ \barh[\sX,\tau,\dB(\fv,r)]
\quad\geq\quad 
\limsup_{r\goto\oo}\ \barh[\sX,\tau,\dB(\fv,r)]
\quad\geeeq{(\dagger)}\quad 
\limsup_{r\goto\oo} \
\frac{S \eps}{4 r} \, \log_2\lb|\sX_{\dB(\fv,r)}\rb|
\\
&=&
\limsup_{r\goto\oo}\
\lb(\frac{S \eps}{4} \cdot
\frac{\log_2\lb|\sX_{\dB(\fv,r)}\rb|}{\lb|\dB(\fv,r)\rb|} \cdot
\frac{\lb|\dB(\fv,r)\rb|}{r}\rb) 
\\&\geq &
\frac{S \eps}{4}  \cdot
\lb(\limsup_{r\goto\oo}\ \frac{\log_2\lb|\sX_{\dB(\fv,r)}\rb|}{\lb|\dB(\fv,r)\rb|}\rb)
\cdot\lb(\liminf_{r\goto\oo}\  \frac{\lb|\dB(\fv,r)\rb|}{r} \rb)
\\&\eeequals{(\ref{entropy.defn})}&
\lb(\frac{S \eps \, \barh_\fv(\sX)}{4} \rb) \
\liminf_{r\goto\oo} \ \frac{\lb|\dB(\fv,r)\rb|}{r} 
\quad\eeequals{(*)}\quad \oo,
\eeq
as desired. 
Here, $(\dagger)$ is by Claim \ref{infinite.tau.entropy.C1}, and
$(*)$ is because $\barh_\fv(\sX)>0$ and
$(\dV,\connect)$ has superlinear connectivity at $\fv$.
\ethmprf

\Lemma{\label{infinite.entropy.means.not.posexpansive}}
{
Let $(\AV,\sX,\Phi)$ symbolic dynamical system with a subsymmetry
$\tau:\dV\into\dV$. If $\barh(\sX,\tau)=\oo$,  then $(\sX,\Phi)$ is not posexpansive.
}
\bthmprf
(by contradiction) \   Suppose
$(\sX,\Phi)$ is posexpansive. Let $\dW_0\subseteq\dV $ be a posexpansive
window, with panorama $\{\dW_0^t\}_{t=0}^\oo$.
For any $n\in\Natur$, let $\dW_n:=\tau^n(\dW_0)$, and 
for all $t\in\Natur$, let
$\dW_n^t := \tau^n(\dW_0^t)$.

\Claim{For all $n\in\Natur$, $\{\dW_n^t\}_{t=0}^\oo$ is the panorama of $\dW_n$.}
\bclaimprf
For any $\bx,\bx'\in\sX$, and any $T\in\Natur$, we have
\beq
\lefteqn{\statement{$\Phi_{\dW_n}^\CC{0...T}(\bx) = \Phi_{\dW_n}^\CC{0...T}(\bx')$}
\ \iff \ 
\statement{$\Phi^t(\bx)_{\dW_n} = \Phi^t(\bx')_{\dW_n}$, \ $\forall \ t\in\CC{0...T}$} }
\\ &\iff&
\statement{$\tau^n\circ\Phi^t(\bx)_{\dW_0} = \tau^n\circ\Phi^t(\bx')_{\dW_0}$,\  $\forall \ t\in\CC{0...T}$}
\\ &\iiiff{(\dagger)}&
\statement{$\Phi^t\circ\tau^n(\bx)_{\dW_0} = \Phi^t\circ\tau^n(\bx')_{\dW_0}$,\  $\forall \ t\in\CC{0...T}$}
\\ &\iff&
 \statement{$\Phi_{\dW_0}^\CC{0...T}[\tau^n(\bx)] = \Phi_{\dW_0}^\CC{0...T}[\tau^n(\bx')]$}
\\ &\iiimplies{(*)}&
\statement{$\tau^n(\bx)_{\dW_0^T} = \tau^n(\bx')_{\dW_0^T}$}
\ \iff\ 
\statement{$\bx_{\dW_n^T} = \bx'_{\dW_n^T}$},
\eeq
as desired.
Here $(*)$ is by statement (\ref{expansive.panorama1}),
because  $\{\dW_0^t\}_{t=0}^\oo$ is the
panorama of $\dW_0$, and $(\dagger)$ is because $\Phi\circ\tau^n=\tau^n\circ\Phi$,
because $\tau$ is a subsymmetry of $(\AV,\sX,\Phi)$.
\eclaimprf
\ignore{
(If $\tau$ is bijective, then $\tau^n$ is an automorphism of $(\sX,\Phi)$;
in this case, $\dW_n$ is itself a posexpansive window for  $(\sX,\Phi)$.  But
we do not require this.)}

\Claim{{\rm(a)} \ There exists $T\in\Natur$
such that $\dW_1\subseteq\dW_0^T$.  

{\rm(b)} \ For all $t\in\Natur$, we have  $\dW_1^t\subseteq\dW_0^{T+t}$.

{\rm(c)}  \ For all $n\in\Natur$, we have $\dW_n\subseteq\dW_0^{nT}$.

{\rm(d)} \  For all $n,t\in\Natur$, we have $\dW^t_n \subseteq \dW_0^{nT+t}$.
}
\bclaimprf (a) follows from equation (\ref{expansive.panorama2}).
To see (b), let $\bx,\bx'\in\sX$.  Then
\beq
\statement{$\Phi_{\dW_0}^\CC{0\ldots T\!+\!t}(\bx) = \Phi_{\dW_0}^\CC{0\ldots T\!+\!t}(\bx')$}
&\iff&
\statement{$\Phi_{\dW_0}^\CC{s\ldots T\!+\!s}(\bx) = \Phi_{\dW_0}^\CC{s\ldots T\!+\!s}(\bx')$  for all $s\in\CC{0...t}$}
\\&\iiiff{(\dagger)}&
\statement{$\Phi_{\dW_0}^\CC{0\ldots T}[\Phi^s(\bx)] = \Phi_{\dW_0}^\CC{0\ldots T}[\Phi^s(\bx')]$ for all $s\in\CC{0...t}$}
\\&\iiimplies{(*)}&
\statement{$\Phi^s(\bx)_{\dW_0^T}=\Phi^s(\bx')_{\dW_0^T}$ for all $s\in\CC{0...t}$}
\\&\iiimplies{(@)}&
\statement{$\Phi^s(\bx)_{\dW_1}=\Phi^s(\bx')_{\dW_1}$ for all $s\in\CC{0...t}$}
\\&\iiiff{(\diamond)}&
\statement{$\Phi_{\dW_1}^\CC{0...t}(\bx)=\Phi_{\dW_1}^\CC{0...t}(\bx')$}
\iiimplies{(\ddagger)}
\statement{$\bx_{\dW_1^t} = \bx'_{\dW_1^t}$}.
\eeq
Thus, $\dW_1^t\subseteq\dW_0^{T+t}$, as desired.  Here,
$(*)$ is by statement (\ref{expansive.panorama1}),
because  $\{\dW_0^t\}_{t=0}^\oo$ is the
panorama of $\dW_0$.\quad
$(\dagger)$  is because $\Phi_{\dW_0}^\CC{s\ldots T\!+\!s}(\bx)=\Phi_{\dW_0}^\CC{0\ldots T}[\Phi^s(\bx)]$ (and likewise for $\bx'$).\quad
$(@)$ is by part (a). Finally,
$(\diamond)$ is by definition of $\Phi_{\dW_1}^\CC{0...t}(\bx)$,
and $(\ddagger)$ is by statement (\ref{expansive.panorama1}) and Claim 1.

(c) (by induction on $n$) \ The case ($n=1$) is Part (a).  Now
suppose inductively that $\dW_n\subseteq \dW_0^{nT}$.
Then $\dW_{n+1} = \tau(\dW_n) \ \subseteeeq{(\dagger)} \ \tau(\dW_0^{nT})
= \dW_{1}^{nT} \ \subseteeeq{(*)} \ \dW_0^{nT+T}=\dW_0^{(n+1)T}$.
Here, $(\dagger)$ is by the induction hypothesis, and $(*)$ is by part (b).

(d) is the same argument as in part (b), but using part (c)
instead of part (a).
\eclaimprf

Now, for any $H>0$, we can find some finite subset $\dF\subset\dV$ such
that $\barh(\sX,\tau,\dF)>H$ (because $\barh(\sX,\tau)=\oo$).
Equation (\ref{expansive.panorama2}) yields some
$t$ such that $\dF\subseteq\dW_0^{t}$.  Thus, for any $N\in\Natur$, and
all $n\in\CC{0...N}$, we have
\[
\tau^n[\dF] \quad \subseteq\quad  \tau^n[\dW_0^{t}]
\quad = \quad 
\dW_n^{t}\quad \subseteeeq{(*)}\quad \dW_0^{nT+t}
\quad \subseteeeq{(\dagger)}\quad \dW_0^{NT+t},
\]
 where $(*)$ is by Claim 2(d), and $(\dagger)$ is because $nT+t\leq NT+t$.
Thus, 
\beqn
\label{infinite.entropy.means.not.posexpansive.e2}
\mbox{if} \quad \dF(N) \quad := \quad \Union_{n=0}^N \tau^n(\dF),
\quad \mbox{then} \quad \dW_0^{NT+t} \ \ \supseteq  \ \  \dF(N).
\eeqn
Thus,
\beqn
\label{infinite.entropy.means.not.posexpansive.e3}
\log_2\lb|\Phi_{\dW_0}^{\CC{0\ldots NT\!+\!t}}(\sX)\rb|
\quad \geeeq{(*)} \quad
\log_2\lb|\sX_{\dW_0^{NT+t}}\rb|
\quad\geeeq{(\dagger)}\quad
\log_2\lb|\sX_{\dF(N)}\rb|.
\eeqn
Here, $(*)$ is because statement (\ref{expansive.panorama1}) yields an
injection from $\sX_{\dW_0^{NT+t}}$ into $\Phi_{\dW_0}^{\CC{0\ldots NT\!+\!t}}(\sX)$;
meanwhile $(\dagger)$ is by equation (\ref{infinite.entropy.means.not.posexpansive.e2}).
Let $\sB:=\sA^{\dW_0}$ and  $\sY:=\Phi_{\dW_0}^\Natur(\sX)\subseteq\BN$.
Then
\begin{eqnarray}
\nonumber
\barh(\sY,\sigma)
&:=&
\limsup_{M\goto\oo} \ \frac{\log_2|\sY_{\CC{0...M}}|}{M}
\quad\geq\quad
\limsup_{N\goto\oo} \ \frac{\log_2|\sY_{\CC{0\ldots NT\!+\!t}}|}{NT+t}
\\&=&
\limsup_{N\goto\oo}\  \frac{\log_2\lb|\Phi_{\dW_0}^\CC{0\ldots NT\!+\!t}(\sX)\rb|}{NT+t}
\quad\geeeq{(\dagger)}\quad
\limsup_{N\goto\oo} \ \frac{\log_2\lb|\sX_{\dF(N)}\rb|}{NT+t} \nonumber
\\ &=&\limsup_{N\goto\oo} \ \lb(\frac{N}{NT+t}\rb)\lb(\frac{\log_2\lb|\sX_{\dF(N)}\rb|}{N}\rb)
\quad\eeequals{(\ddagger)}\quad
\frac{\barh(\sX,\tau,\dF)}{T}
\quad>\quad \frac{H}{T}.\qquad
\label{infinite.entropy.means.not.posexpansive.e4}
\end{eqnarray}
Here, 
$(\dagger)$ is by equation (\ref{infinite.entropy.means.not.posexpansive.e3}),
and $(\ddagger)$ is by defining equation (\ref{tau.entropy.defn1}).
Now, $H$ can be made arbitrarily large, because $\barh(\sX,\tau)=\oo$.
Thus, letting $H\goto\oo$ in equation 
(\ref{infinite.entropy.means.not.posexpansive.e4}),
we conclude that $\barh(\sY,\sigma)=\oo$.

But clearly, $\barh(\sY,\sigma) \leq \log_2|\sB| = \log_2|\sA^{\dW_0}|
= |\dW_0|\cdot\log_2|\sA|<\oo$, because $\sA$ and $\dW_0$ are finite.  Contradiction.
\ethmprf

\bthmprf[Proof of Theorem \ref{not.posexpansive}.]
Combine Lemmas \ref{infinite.tau.entropy} and \ref{infinite.entropy.means.not.posexpansive}.
\ethmprf

\paragraph{Remarks.}
(a) Observe that Lemma \ref{infinite.entropy.means.not.posexpansive}
is really a statement about $(\sX,\Phi)$ as an abstract Cantor
dynamical system with a subsymmetry; it does not depend on any
specific representation of $(\sX,\Phi)$ as a symbolic dynamical
system (i.e. any specific embedding $\sX\subset\AV$ for some $\sA$ and $\dV$).
As such, Lemma \ref{infinite.entropy.means.not.posexpansive}
is an interesting result in itself.

(b) Theorem \ref{not.posexpansive} applies even if $\tau$ and its
iterates are the only symmetries of $(\AV,\sX,\Phi)$.  In particular, we
do {\em not} require the symmetry group of $(\AV,\sX,\Phi)$ to itself
have growth
dimension greater than $1$.

(c) The `weak independence' condition in Theorem \ref{not.posexpansive}
and Lemma \ref{infinite.tau.entropy} is probably not necessary.\hfill$\diamondsuit$.

\section{Propagation, Sensitivity, and Equicontinuity \label{S:propagation}}

If $\Phi:\AV\into\AV$ is continuous, then for all $t\in\Natur$, the
function $\Phi^t:\AV\into\AV$ is also continuous; hence we can apply
Lemma \ref{nhood.corollary} to define input neighbourhoods
$\Phi^t_\In(\fv)\subset\dV$ for all $\fv\in\dV$.  The {\dfn
propagation} of $\Phi$ at $\fv$ is the function
$\rho_\fv:\Natur\into\Natur$ defined by
\beqn
\rho_\fv(T):=
\lb|\Phi^\CC{0...T}_\In(\fv)\rb|,
\quad\mbox{where}\quad
\Phi^\CC{0...T}_\In(\fv) \ := \ \Union_{t=0}^T \Phi^t_\In(\fv),
\quad\mbox{for all $T\in\Natur$.} 
\label{propagation.defn}
\eeqn
Clearly, $\rho_\fv(T)\leq |\dB(\fv,T)|$ (because for all $t\leq T$, we have
$\Phi^t_\In(\fv)\ \subseteq\ \dB(\fv,T)$).  In general, this
inequality may be strict.

Let $(\AV,\sX,\Phi)$ be a symbolic dynamical system,
and let $\fv\in\dV$.  A point $\bx\in\sX$
is {\dfn $\fv$-sensitive}
if there exists a sequence $\{\bx^n\}_{n=1}^\oo\subset\sX$
with $\D \lim_{n\goto\oo}\ \bx^n \ = \ \bx$,  such that
\beqn
\label{sensitivity}
\mbox{For all $n\in\Natur$, there is some $t\in\Natur$ with
 $\Phi^t(\bx^n)_\fv\neq \Phi^t(\bx)_\fv$.}
\eeqn 
We say that $\bx$ is a {\dfn sensitive point} if it is $\fv$-sensitive
for some $\fv\in\dV$.
(If $d$ is any compatible metric on $\sX$, then there is some $\eps$ such that,
for all $\bx,\by\in\sX$, \ $\statement{$x_\fv\neq y_\fv$}\implies 
\statement{$d(\bx,\by)>\eps$}$; thus, this definition
is equivalent to the ordinary metric definition of `sensitivity').

\Proposition{\label{sensitivity.vs.propagation1}}
{
Let $(\AV,\sX,\Phi)$ be a symbolic dynamical system.
\bthmlist
\item  Let $\fv\in\dV$ and suppose $\Phi$ has propagation $\rho_\fv$
at $\fv$.  Then
$\statement{$\rho_\fv$ is unbounded}
\iff
\statement{$(\sX,\Phi)$ 
has a $\fv$-sensitive point}$.

\item $(\sX,\Phi)$ has a sensitive point
if and only if
there exists $\fv\in\dV$ with
unbounded propagation.
\ethmlist
}
\bthmprf (a)
``$\implies$'' \quad
For all $r\in\Natur$, there exists $T(r)\in\Natur$ such that
$\rho_\fv[T(r)]>|\dB(\fv,r)|$, which means there is some $t(r)\in\CC{0...T(r)}$ with
$\Phi^{t(r)}_\In(\fv)\not\subseteq\dB(\fv,r)$.  Let 
$\fw\in\Phi^{t(r)}_\In(\fv)\setminus\dB(\fv,r)$.
The local rule $\phi^{t(r)}_\fv:\sA^{\Phi^{t(r)}_\In(\fv)}\into\sA$ is proper, so
there exist $\by^r,\bz^r\in\sX$ such that 
\[
{\bf(a)}
\quad y^r_\fu\ = \ z^r_\fu, \ \mbox{for all $\fu\in\dV\setminus\{\fw\}$; \quad but}
 \quad {\bf(b)} \quad \Phi^{t(r)}(\by^r)_\fv \ \neq \  \Phi^{t(r)}(\bz^r)_\fv.   
\]
 Now, $\dB(\fv,r)\subseteq\dV\setminus\{\fw\}$ by construction, so condition (a) means that 
\beqn
\label{fooz}
\by^r_{\dB(\fv,r)}\quad=\quad\bz^r_{\dB(\fv,r)}.
\eeqn
Since $\sX$ is compact, we find some increasing
sequence $\{r_n\}_{n=1}^\oo\in\Natur$ such that the subsequence
$\{\by^{r_n}\}_{n=1}^\oo$ converges in $\sX$ to some point $\bx$.  
Equation (\ref{fooz}) implies that the subsequence $\{\bz^{r_n}\}_{n=1}^\oo$
also converges to $\bx$.  But for all $n\in\Natur$, condition (b) says that
$\Phi^{t(r_n)}(\by^{r_n})_\fv\neq \Phi^{t(r_n)}(\bz^{r_n})_\fv$,
which means that either
(i) $\Phi^{t(r_n)}(\by^{r_n})_\fv\neq \Phi^{t(r_n)}(\bx)_\fv$ or (ii) $\Phi^{t(r_n)}(\bz^{r_n})_\fv\neq \Phi^{t(r_n)}(\bx)_\fv$
(or both).
In case (i), define $\bx^n:=\by^{r_n}$, while in case (ii), define $\bx^n:=\bz^{r_n}$;
then we obtain a sequence $\{\bx^n\}_{n=1}^\oo$ converging to $\bx$,
and satisfying  condition (\ref{sensitivity}).

``$\seilpmi$'' \ For any $R\in\Natur$, we must find some $T\in\Natur$
such that $\rho_\fv(T)>R$.
Let $\bx\in\sX$ be a $\fv$-sensitive point, and let  $\{\bx^n\}_{n=1}^\oo\subset\sX$ be a sequence converging to $\bx$
and satisfying  condition (\ref{sensitivity}).  
Now, there exist some $n\in\Natur$ such that $\bx^n_{\dB(\fv,R)}=\bx_{\dB(\fv,R)}$;
but there also exists some $T$ such that $\Phi^T(\bx^n)_\fv\neq \Phi^T(\bx)_\fv$.
This means that $\Phi^T_\In(\fv)\not\subseteq \dB(\fv,R)$.  Let $\fw_T\in 
\Phi^T_\In(\fv)\setminus \dB(\fv,R)$.  Then there exists a
directed path $\fw_T\connect\fw_{T-1}\connect\cdots\connect \fw_2\connect\fw_1\connect\fv$
such that $\fw_t\in \Phi^t_\In(\fv)$ for all $t\in\CC{1...T}$.  Furthermore, this path
must have length $L>R$ (even after removing repeated entries), because $\fw_T\not\in\dB(\fv,R)$.  
Thus, $\rho_\fv(T) \ := \ 
\lb|\Phi^\CC{0..T}_\In(\fv)\rb|
\ \geq\  \lb|\{\fw_T,\ldots,\fw_2,\fw_1\}\rb| \ =\  L\ >\  R$,
as desired.  This works for any $R\in\Natur$; hence the function $\rho_\fv$ is unbounded.

(b) follows immediately from (a).
\ethmprf

Let $\dW\subset\dV$ be some finite subset.
We say that $\Phi$ is {\dfn $\dW$-equicontinuous} 
if there exists a finite subset $\dU\subset\dV$
containing $\dW$ (called the {\dfn envelope} of $\dW$), such that:
\beqn
\label{W.equicontinuity}
\mbox{For all $\bx,\by\in\sX$,}\qquad
\statement{$\by_\dU=\bx_\dU$}\quad \implies\quad
\statement{$\Phi^t(\by)_\dW =\Phi^t(\bx)_\dW, \ \forall t\in\Natur$}.
\eeqn 
We say that $\Phi$ is {\dfn equicontinuous} if $\Phi$ is 
$\dW$-equicontinuous for every finite subset $\dW\subset\dV$.
(If $d$ is any compatible metric on $\sX$, then for any $\eps>0$ there
is some finite subset $\dW\subset\dV$ such that
for all $\bx,\by\in\sX$, \ $\statement{$\bx_\dW= \by_\dW$}\implies 
\statement{$d(\bx,\by)<\eps$}$.  Likewise, for any finite $\dU\subset\dV$,
there is some $\del>0$ such that 
for all $\bx,\by\in\sX$, \ $\statement{$d(\bx,\by)<\del$}
\implies \statement{$\bx_\dU= \by_\dU$}$.  Thus, our definition
is equivalent to the ordinary metric definition of `equicontinuity').

A topological dynamical system $(\sX,\Phi)$ is an {\dfn odometer}
if $(\sX,\Phi)$ is an inverse limit of a sequence of finite, periodic
dynamical systems.  That is:
\beqn
\label{inverse.limit}
(\sX,\Phi)
\quad:=\quad
\lim\lb\{ \cdots \stackrel{\pi_3}{\longrightarrow} (\sX_3,\phi_3)
\stackrel{\pi_2}{\longrightarrow} (\sX_2,\phi_2)
\stackrel{\pi_1}{\longrightarrow} (\sX_1,\phi_1) \rb\},
\eeqn
where, for all $n\in\Natur$, $\sX_n$ is a finite set, $\phi_n:\sX_n\into\sX_n$
is a cyclic permutation, and $\pi_n:(\sX_{n+1},\phi_{n+1})\into(\sX_n,\phi_n)$ is a factor mapping.

For any $x\in\sX$, let $\barsO_x:=\overline{\set{\Phi^t(x)}{t\in\Natur}}$
be the {\dfn $\Phi$-orbit closure} of $x$; then $(\barsO_x,\Phi)$
is itself a topological dynamical system.
The system $(\sX,\Phi)$ is an {\dfn odometer bundle}
if, for every $x\in\sX$, the system $(\barsO_x,\Phi)$ is
an odometer.  Thus, $(\sX,\Phi)$ can be decomposed into a 
(possibly infinite) disjoint union of (possible non-isomorphic) odometers.

For example, for all $n\in\Natur$, let $\sX_n$ be a finite set, and
let $\phi_n:\sX_n\into\sX_n$ be a permutation (possibly
with multiple disjoint orbits).  Suppose $(\sX,\Phi)$ arises as
the inverse limit (\ref{inverse.limit});  then $(\sX,\Phi)$ is
an odometer bundle. 

\Proposition{\label{equicontinuity.vs.propagation}}
{
For all $\fv\in\dV$, let $\rho_\fv:\Natur\into\Natur$ be the propagation
of $\Phi$ at $\fv$.
\bthmlist  
  \item Let $\dW\subset\dV$ be a finite subset.  If $\rho_\fw$ is
bounded for all $\fw\in\dW$, then $\Phi$ is $\dW$-equicontinuous.

  \item $\statement{$\rho_\fv$ is
bounded for all $\fv\in\dV$}\iff \statement{$\Phi$ is equicontinuous}$.

  \item  If $\Phi:\sX\into\sX$ is equicontinuous and surjective, then $(\sX,\Phi)$ is
an odometer bundle.
\ethmlist
}
\bthmprf
(a) For all $\fw\in\dW$, there is some $R_\fw$ such that
$\rho_\fw(t)<R_\fw$ for all $t\in\Natur$.  Let $\D R:=\max_{\fw\in\dW}
R_\fw$; then $R$ is finite because $\dW$ is finite.  Let 
$\dU:=\dB(\dW,R)$.  

\Claim{For all $\fw\in\dW$, and all $t\in\Natur$, we have $\Phi^t_\In(\fw)\subseteq\dU$.}
\bclaimprf
(by contradiction) \  Suppose $\Phi^t_\In(\fw)\not\subseteq\dU$.
Let $\fv\in\Phi^t_\In(\fw)\setminus\dU$,
and just as in the proof of Proposition \ref{sensitivity.vs.propagation1} ``$\seilpmi$'',
construct a path from $\fv$ to $\fw$ of length $L>R$.
Conclude that $\rho_\fw(t) >R$.   Contradiction.
\eclaimprf

Suppose $\bx_\dU=\by_\dU$.  Then
for all $\fw\in\dW$, and all $t\in\Natur$, \ Claim 1 implies that
$\Phi^t(\bx)_\fw=\Phi^t(\bx)_\fw$.  In other words, 
$\Phi^t(\by)_\dW =\Phi^t(\bx)_\dW$,  for all $t\in\Natur$.
Thus, $\dU$ is an envelope for $\dW$.

 (b) \ ``$\implies$'' follows immediately from part (a).   For ``$\seilpmi$'',
note that an equicontinuous system can have no sensitive points; now
apply the contrapositive of Proposition \ref{sensitivity.vs.propagation1}(b).

(c) \  Let $\dW_1\subset\dW_2\subset\dW_3\subset\cdots$
be an increasing sequence of finite sets, with $\Union_{n=1}^\oo \dW_n=\dV$.
For all $n\in\Natur$, let $\sB_n:=\sX_{\dW_n}$, 
and let $\Phi^\Natur_{\dW_n}:\sX\into\sB_n^\Natur$ be as in \S\ref{S:posexp}.
Let $\sY_n:=\Phi^\Natur_{\dW_n}(\sX)\subseteq\sB_n^\Natur$.  
If $\sigma_n:\sB_n^\Natur\into\sB_n^\Natur$ is the
shift map, then $\sigma_n\circ \Phi^\Natur_{\dW_n} = \Phi^\Natur_{\dW_n}\circ \Phi$.
Furthermore, $\sigma_n(\sY_n)=\sY_n$, because $\Phi$ is surjective.

For all $n\in\Natur$, let $\tlpi_n:\sX_{\dW_{n+1}}\into\sX_{\dW_n}$
be the projection (i.e. $\tlpi_n(\bx_{\dW_{n+1}}):=\bx_{\dW_n}$
for all $\bx\in\sX$).  Define
$\pi_n:\sY_{n+1}\into\sY_n$ as follows:
if $\by\in\sY_{n+1}$, write $\by= [\by^t]_{t=0}^\oo$ 
where $\by^t\in\sX_{\dW_{n+1}}$ for all $t\in\Natur$;
then define $\pi_n(\by):=[\tlpi_n(\by^t)]_{t=0}^\oo$.  Clearly,
$\pi_n:(\sY_{n+1},\sigma_{n+1})\into(\sY_n,\sigma_n)$ is a factor
mapping, and 
$(\sX,\Phi)
\ = \ 
\lim\lb\{ \cdots \stackrel{\pi_3}{\longrightarrow} (\sY_3,\sigma_3)
\stackrel{\pi_2}{\longrightarrow} (\sY_2,\sigma_2)
\stackrel{\pi_1}{\longrightarrow} (\sY_1,\sigma_1) \rb\}$.

Everything so far is true for any symbolic dynamical
system.  Now we use equicontinuity.

\claim{For all $n\in\Natur$,  \ $\sY_n$ is finite and $\sigma_n:\sY_n\into\sY_n$ is a permutation.}
\bclaimprf
Let $\dU_n\subset\dV$ be the envelope of $\dW_n$ (a finite set).
For any $\bx,\bx'\in\sX$, if $\bx_{\dU_n}=\bx'_{\dU_n}$,
then $\Phi^\Natur_{\dW_n}(\bx)=\Phi^\Natur_{\dW_n}(\bx')$.
Thus, $|\Phi^\Natur_{\dW_n}(\sX)|\leq |\sX_{\dU_n}|$
---in other words,  $|\sY_n|\leq |\sX_{\dU_n}|$.
But $|\sX_{\dU_n}|$ is finite because $\dU_n$ is finite.
Thus, $\sY_n$ is finite.  Thus, $\sigma_n$ is bijective (because
it is surjective).
\eclaimprf
Thus, we have represented $(\sX,\Phi)$ as an inverse limit of
finite permutation dynamical systems;  thus, $(\sX,\Phi)$ is
an odometer bundle.
\ethmprf

For example: any symbolic dynamical system with the network  in
Figure \ref{fig:odometer.network} must be equicontinuous.

\ignore{\Remark{
 If $(\sX,d)$ is a Cantor space with a compatible metric,
and  $\Phi:\sX\into\sX$ is surjective,
then Propositions \ref{sensitivity.vs.propagation1} and \ref{equicontinuity.vs.propagation}(b) imply that
either  $(\sX,\Phi)$ has a sensitive point, or $(\sX,\Phi)$ is conjugate
to an odometer bundle.}}

\section{An expansive system of dimension two
\label{S:counterexample}}

  The symmetry condition in Theorem \ref{not.posexpansive} is probably
not necessary.  However, some sort of condition is required beyond
merely superlinear connectivity.  To demonstrate this, we will
construct an example of a symbolic dynamical system which is
posexpansive, despite having connectivity dimension two.

\begin{figure}[t]
\centerline{
\psfrag{1}[][]{\tiny 0}
\psfrag{2}[][]{\tiny 1}
\psfrag{3}[][]{\tiny 2}
\psfrag{4}[][]{\tiny 3}
\psfrag{5}[][]{\tiny 4}
\psfrag{6}[][]{\tiny 5}
\psfrag{7}[][]{\tiny 6}
\psfrag{8}[][]{\tiny 7}
\psfrag{9}[][]{\tiny 8}
\psfrag{10}[][]{\tiny 9}
\psfrag{11}[][]{\tiny 10}
\psfrag{12}[][]{\tiny 11}
\psfrag{13}[][]{\tiny 12}
\psfrag{14}[][]{\tiny 13}
\psfrag{15}[][]{\tiny 14}
\psfrag{16}[][]{\tiny 15}
\psfrag{17}[][]{\tiny 16}
\psfrag{18}[][]{\tiny 17}
\psfrag{19}[][]{\tiny 18}
\psfrag{20}[][]{\tiny 19}
\psfrag{21}[][]{\tiny 20}
\psfrag{22}[][]{\tiny 21}
\psfrag{23}[][]{\tiny 22}
\psfrag{24}[][]{\tiny 23}
\psfrag{25}[][]{\tiny 24}
\psfrag{26}[][]{\tiny 25}
\psfrag{27}[][]{\tiny 26}
\psfrag{28}[][]{\tiny 27}
\psfrag{29}[][]{\tiny 28}
\psfrag{30}[][]{\tiny 29}
\psfrag{31}[][]{\tiny 30}
\psfrag{32}[][]{\tiny 31}
\psfrag{33}[][]{\tiny 32}
\psfrag{34}[][]{\tiny 33}
\psfrag{35}[][]{\tiny 34}
\psfrag{36}[][]{\tiny 35}
\includegraphics[scale=0.4,angle=-90]{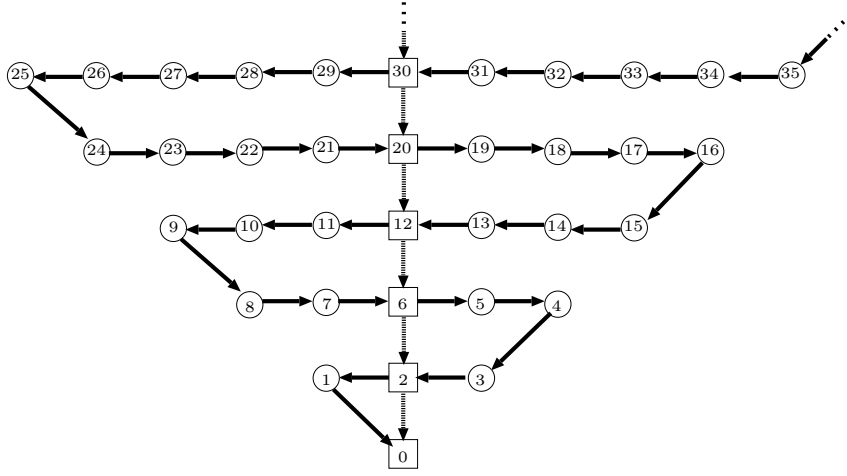}
}
\caption{\label{fig:expansive.2d} An expansive, two-dimensional
symbolic dynamical system.}
\end{figure}

Let $\sA:=\Zahlmod{2}\x\Zahlmod{2}$, and let $\dV$ be the digraph
shown in Figure \ref{fig:expansive.2d}.  Let $\dV_\Box$ be the set of
vertices indicated by boxes and indexed by
$\dM:=\{0,2,6,12,20,\ldots,m_k,\ldots\}$, where $m_k :=
\sum_{i=0}^{k} 2j$.  We denote the $m_k$th square vertex by
$\Box_{m_k}$.  Let $\dV_\circ$ be the set of vertices indicated by
circles; then $\dV=\dV_\Box\disj\dV_\circ$.  We denote the $n$th
element of $\dV_\circ$ by $\bigcirc_n$.  

For any vertex $\fv_n\in\dV$,
the state of $\fv_n$ is an
ordered pair $\lb({a_n \atop b_n}\rb)$, where $a_n,b_n\in\Zahlmod{2}$.
Let $\sX:=\set{\ba\in\AV}{b_n=0, \ \forall \ \fv_n\in\dV_\circ}$.
Thus, if $\sX_\fv$ is the projection of $\sX$ onto vertex $\fv$, then
we have $\sX_\fv=\Zahlmod{2}\x\Zahlmod{2}$ if $\fv\in\dV_\Box$,
and $\sX_\fv=\Zahlmod{2}\x\{0\}$ if $\fv\in\dV_\circ$.
(Clearly $h(\sX)\geq \log_2(2)=1$, and $\sX$ has weak independence.)

 The local rule of each cell depends entirely upon its one or two
input cells, and not on itself, as follows. 
For any $\bigcirc_n\in\dV_\circ$,
we define
$\phi_n:\sX_{n+1}\into\sX_n$ by $\phi_n\lb({a_{n+1}\atop b_{n+1}}\rb) \ =
 \lb({a_{n+1}\atop 0}\rb)$ ---that is, $\phi_n$ simply copies the
first coordinate of $\bigcirc_{n+1}$ (or $\Box_{n+1}$) into  $\bigcirc_n$.

  The cell $\Box_{m_k}$ is connected to both  
$\bigcirc_{(m_k)+1}$ and  $\Box_{m_{(k+1)}}$.  Its local rule $\phi_{m_k}:\sX_{(m_k)+1}
\x \sX_{m_{(k+1)}}\into\sX_{m_k}$ is defined as follows:
\beqn
\label{local.rule.for.box.cell}
  \phi_{m_k} \lb( \lb({a_{(m_k)+1}\atop 0}\rb) \ , \  \lb({a_{m_{(k+1)}} \atop b_{m_{(k+1)}}}\rb)\rb)
\quad:=\quad
\lb({a_{(m_k)+1} \atop {a_{m_{(k+1)}} + b_{m_{(k+1)}} }}\rb)
\eeqn

\Lemma{\label{posexpansive.counterexample1}}
{The system $(\sX,\Phi)$ is posexpansive\footnote{See the start of \S\ref{S:posexp} for the definitions of `posexpansive' and
`window'.}, with posexpansive window $\{\Box_{0}\}$.}
\bthmprf
Define $\Gam_1:\sX\into\Zahlmod{2}^\Natur$ and $\Gam_2:\sX\into\Zahlmod{2}^\dM$ by
\beq
\Gam_1\lb[\lb({a_0\atop b_0}\rb), \ \lb({a_1\atop 0}\rb), \ \lb({a_2\atop b_2}\rb), \ \lb({a_3\atop 0}\rb), \ \lb({a_4\atop 0}\rb),  \ \lb({a_5\atop 0}\rb),   \ \lb({a_6\atop b_6}\rb), \ \ldots\rb]
&:=&
[a_0, a_1, a_2, a_3, a_4, a_5, \ldots];
\\
\Gam_2\lb[\lb({a_0\atop b_0}\rb), \ \lb({a_1\atop 0}\rb), \ \lb({a_2\atop b_2}\rb), \ \lb({a_3\atop 0}\rb), \ \lb({a_4\atop 0}\rb),  \ \lb({a_5\atop 0}\rb),   \ \lb({a_6\atop b_6}\rb), \ \ldots\rb]
&:=&
[b_0, b_2, b_6, b_{12}, b_{20}, \ldots].
\eeq
Thus, $\Gam:=\Gam_1\x\Gam_2:\sX\into\Zahlmod{2}^\Natur\x\Zahlmod{2}^\dM$ is a bijection.
Let $\tl\Phi:=\Gam\circ\Phi\circ\Gam^{-1}:\Zahlmod{2}^\Natur\x\Zahlmod{2}^\dM\into\Zahlmod{2}^\Natur\x\Zahlmod{2}^\dM$;
then $\Phi$ is conjugate to $\tlPhi$ via $\Gam$, so it
suffices to show that the system $(\Zahlmod{2}^\Natur\x\Zahlmod{2}^\dM,\tlPhi)$ is posexpansive.

\ignore{Let $\sigma:\Zahlmod{2}^\Natur\into\Zahlmod{2}^\Natur$ be the unilateral shift.  Observe that
\beqn
\label{shift.on.first.coordinate}
\Gam_1\circ\Phi \quad=\quad \sigma\circ\Gam_1.
\eeqn}
Let $(\ba^0,\bb^0)\in \Zahlmod{2}^\Natur\x\Zahlmod{2}^\dM$ be some initial state,
and let $(\ba^t,\bb^t):=\tlPhi^t(\ba^0,\bb^0)$ for all $t\in\Natur$.
Let $J\in\Natur$, and consider the sequence of
the first $m_J$ states of cell $\Box_0$:
\[
\bx_J\quad:=\quad
\lb( {a^0_0 \atop b^0_0}, \ {a^1_0 \atop b^1_0}, \ \ldots,
{a^{m_J}_0 \atop b^{m_J}_0} \rb).
\]
Clearly, the information in $\bx_J$ is sufficient to determine
the values of $a^0_0, a^0_1, a^0_2,\ldots,a^0_{m_J}$, 
because for all $t\in\Natur$ we have $a^t_0 = a^0_{t}$,
because $\tlPhi$ just acts like a unilateral shift
on the `$a$' coordinates.
It remains to show that $\bx_J$ is also sufficient to determine
the values of $b^0_0, b^0_2, b^0_6,\ldots,b^0_{m_J}$. 

\Claim{For any $j\in\CC{1...J}$, and any $t\in\CC{0\ldots m_J\!-\!m_j}$, 
the information in $\bx_J$ determines the value of $b^t_{m_j}$.}
\bclaimprf
(by induction on $j$) \ 

{\em Base case} ($j=1$). \ Note that $m_1=2$. 
We are given $b^0_0,b^1_0,\ldots,b^{m_J}_0$ in the bottom row of $\bx_J$.
For any $t\in\CC{0\ldots m_J\!-\!2}$, we also know the value of
$a^{t}_2=a^{t+2}_0$, and thus, we can compute
$b^t_2 = b^{t+1}_0 - a^t_2$, because
 $b^{t+1}_0 = a^t_2 + b^t_2$, because substituting $m_0=0$ and $m_1=2$ into
eqn.(\ref{local.rule.for.box.cell}) yields
\[
\lb({a^{t+1}_0 \atop b^{t+1}_0 }\rb)
\quad=\quad
  \phi_{0} \lb( \lb({a^t_{1}\atop0}\rb) \ , \  \lb({a^t_{2} \atop b^t_{2}}\rb)\rb)
\quad:=\quad
\lb({a^t_{1} \atop {a^t_2 + b^t_{2} }}\rb).
\]
{\em Induction.} \  Fix $j\in\CC{1...J}$, and suppose we know the values of 
$b^t_{m_j}$ for all $t\in\CC{0\ldots m_J\!-\!m_{j}}$.
Then we know 
$b^{t+1}_{m_j}$ for all $t\in\CC{0\ldots m_J\!-\!m_{(j+1)}}$
(because $m_J-m_{(j+1)}+1 \ \leq \  m_J-m_{j}$
because $m_{(j+1)} > m_j + 1$).
For any  $t\in\CC{0\ldots m_J\!-\!m_{(j+1)}}$,
we also know $a^{t}_{m_{(j+1)}}=a^{t+{m_{(j+1)}}}_0$
and thus, we can compute
$b^t_{m_{(j+1)}} = b^{t+1}_{m_j} - a^t_{m_{(j+1)}}$, because
$b^{t+1}_{m_j} = a^t_{m_{(j+1)}} + b^t_{m_{(j+1)}}$, because
 eqn.(\ref{local.rule.for.box.cell}) says
\[
\lb({a^{t+1}_{m_j} \atop b^{t+1}_{m_j} }\rb)
\quad=\quad
  \phi_{m_j} \lb( \lb({a^t_{m_j+1}\atop 0}\rb) \ , \  \lb({a^t_{m_{(j+1)}} \atop b^t_{m_{(j+1)}}}\rb)\rb)
\quad:=\quad
\lb({a^t_{m_j+1} \atop {a^t_{m_{(j+1)}} + b^t_{m_{(j+1)}} }}\rb).
\]
\eclaimprf

In particular, Claim 1 implies that,
for any $j\in\CC{1...J}$, 
the information in $\bx_J$ determines $b^0_{m_j}$.
Thus, given $\bx_J$, we can recover
$a^0_0, a^0_2, a^0_3,\ldots,a^0_{m_J}$  and also
$b^0_0, b^0_2, b^0_6,\ldots,b^0_{m_J}$.  This works
for all $J\in\Natur$;  thus, $\{\Box_0\}$ is
a posexpansive window for $(\sX,\Phi)$.
\ethmprf

\begin{figure}[t]
\centerline{
\psfrag{1}[][]{\tiny 0}
\psfrag{2}[][]{\tiny 1}
\psfrag{3}[][]{\tiny 2}
\psfrag{4}[][]{\tiny 3}
\psfrag{5}[][]{\tiny 4}
\psfrag{6}[][]{\tiny 5}
\psfrag{7}[][]{\tiny 6}
\psfrag{8}[][]{\tiny 7}
\psfrag{9}[][]{\tiny 8}
\psfrag{10}[][]{\tiny 9}
\psfrag{11}[][]{\tiny 10}
\psfrag{12}[][]{\tiny 11}
\psfrag{13}[][]{\tiny 12}
\psfrag{14}[][]{\tiny 13}
\psfrag{15}[][]{\tiny 14}
\psfrag{16}[][]{\tiny 15}
\psfrag{17}[][]{\tiny 16}
\psfrag{18}[][]{\tiny 17}
\psfrag{19}[][]{\tiny 18}
\psfrag{20}[][]{\tiny 19}
\psfrag{21}[][]{\tiny 20}
\psfrag{22}[][]{\tiny 21}
\psfrag{23}[][]{\tiny 22}
\psfrag{24}[][]{\tiny 23}
\psfrag{25}[][]{\tiny 24}
\psfrag{26}[][]{\tiny 25}
\psfrag{27}[][]{\tiny 26}
\psfrag{28}[][]{\tiny 27}
\psfrag{29}[][]{\tiny 28}
\psfrag{30}[][]{\tiny 29}
\psfrag{31}[][]{\tiny 30}
\psfrag{32}[][]{\tiny 31}
\psfrag{33}[][]{\tiny 32}
\psfrag{34}[][]{\tiny 33}
\psfrag{35}[][]{\tiny 34}
\psfrag{36}[][]{\tiny 35}
\psfrag{A}[][]{$\Phi^{\{0,1\}}_\In(\Box_0)$}
\psfrag{B}[][]{$\Phi^\CC{0...2}_\In(\Box_0)$}
\psfrag{C}[][]{$\Phi^\CC{0...3}_\In(\Box_0)$}
\psfrag{D}[][]{$\Phi^\CC{0...4}_\In(\Box_0)$}
\psfrag{E}[][]{$\Phi^\CC{0...5}_\In(\Box_0)$}
\psfrag{F}[][]{$\Phi^\CC{0...6}_\In(\Box_0)$}
\psfrag{G}[][]{$\Phi^\CC{0...7}_\In(\Box_0)$}
\includegraphics[scale=0.2]{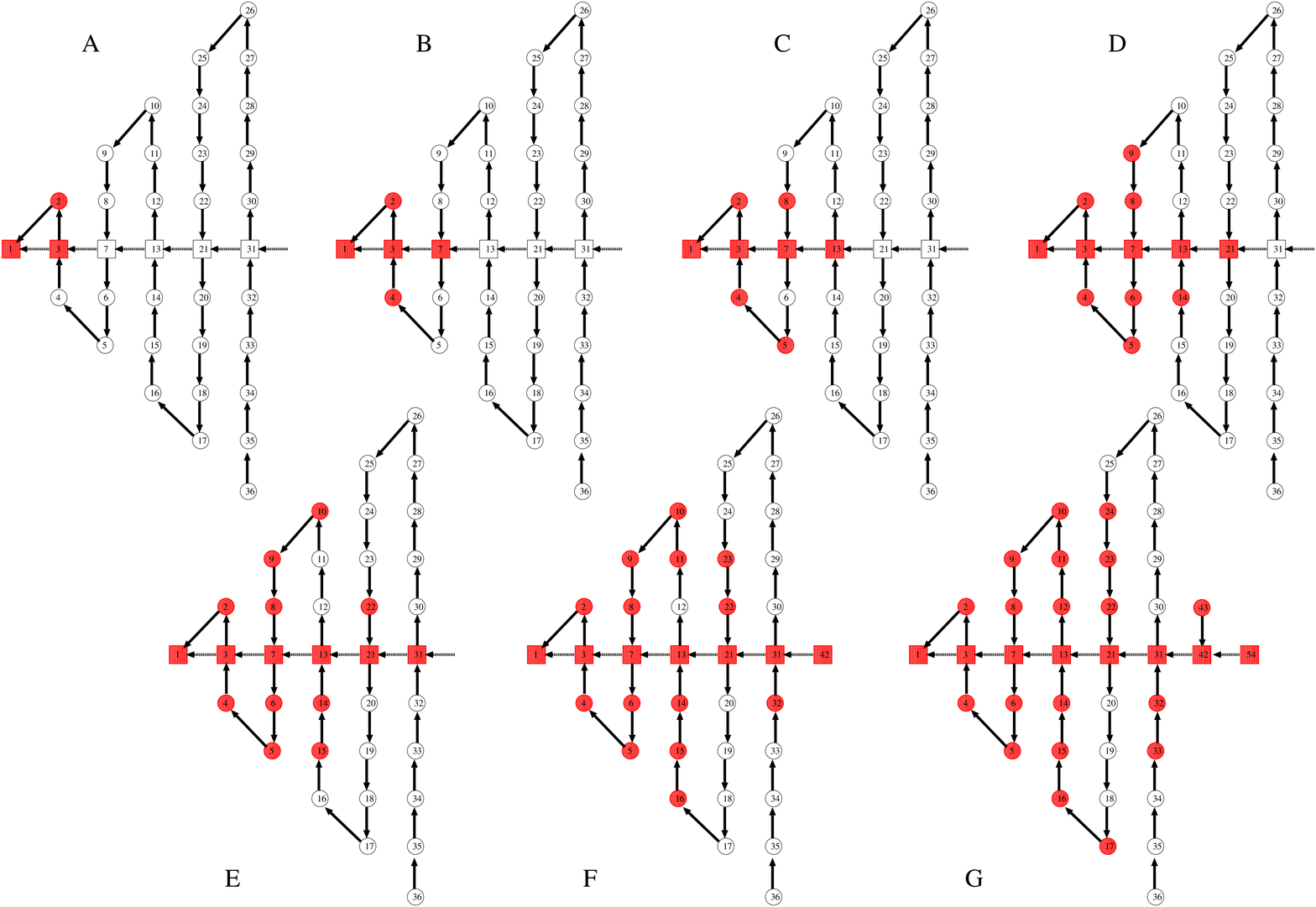}
}
\caption{The quadratic growth of $\Phi^\CC{0...T}_\In(\Box_0)$.
\label{fig:expansive.2d.ball}} 

\end{figure}

\Lemma{\label{posexpansive.counterexample2}}
{
{\rm(a)}\ $\Phi$ has quadratic propagation\footnote{See equation (\ref{propagation.defn}) in  \S\ref{S:propagation} for the definition
of `propagation'.}
 {\rm($\D \liminf_{t\goto\oo}\ (\rho_\fv(t)/t^2) > 0$)} at every $\fv\in\dV_\Box$.

{\rm(b)}\ $\dim_\fv(\dV,\connect)= 2$ for all $\fv\in\dV$.
}
\bthmprf 
(a) We will show this at vertex $\Box_0$;  the proof
at other vertices is similar.
  Observe that $\Phi^\CC{0...T}_\In(\Box_0)$ grows quadratically as
$T\goto\oo$, as shown in Figure \ref{fig:expansive.2d.ball}.
To be precise, for any $T>0$,
 $\Phi^\CC{0...T}_\In(\Box_0)$ contains the cells
$\{\Box_0,\Box_{2},\Box_{6},\ldots,\Box_{m_T}\}$,
and also contains the cells $\bigcirc_{m_t+s}$ for each $t\in\CC{1...T}$ and
$s\in\CC{1\ldots T\!-\!t}$.

(b) From (a) we have  $\dim_\fv(\dV,\connect)\geq 2$ for all $\fv\in\dV_\Box$
(because $|\dB(\fv,r)|\geq \rho_\fv(r)$ for all $r\in\Natur$).
But every vertex is downstream from some element of $\dV_\Box$;
thus, Lemma \ref{constant.nhood.growth.rate} implies that  $\dim_\fv(\dV,\connect)\geq 2$
for all $\fv\in\dV$.  On the other hand, it is easy to
see that  $\dim_\fv(\dV,\connect)\leq 2$.
\ethmprf

\section{Lipschitz metrics on Cantor systems \label{S:lipschitz}}

The counterexample of \S\ref{S:counterexample} shows that connectivity
dimension is {\em not} invariant under topological conjugacy:
some systems with dimension two are conjugate to 
subshifts of $(\sA^\Natur,\sigma)$, and the system
$(\sA^\Natur,\sigma)$ has dimension one.
(Likewise, \S\ref{S:counterexample} shows that the growth
rate of the propagation function $\rho$
is not a conjugacy invariant.)\
 However, we will now show that
connectivity dimension {\em is} invariant under a slightly refined
notion of conjugacy, once we impose a suitable metric structure on the
pattern space $\sX$ (see Corollary
\ref{network.dimension.invariant.under.holder.conjugacy} below).

\breath

Let $\sX$ be any Cantor space, and let $d:\sX^2\into\Real_+$ be
a metric compatible with the Cantor topology on $\sX$. The pair
$(\sX,d)$ will be called a {\dfn Cantor metric space}.
Let $\Phi:\sX\into\sX$ be any continuous self-map. 
We say $\Phi$ is {\dfn $d$-Lipschitz} if there is a constant
$\lam>0$ such that, for any $\bx,\bx'\in\sX$,
\beqn
\label{Lipschitz}
d\lb(\Phi(\bx),\Phi(\bx')\rb) \quad \leq \quad \lam \cdot d(\bx,\bx').
\eeqn
In this case, $d$ is called a {\dfn Lipschitz metric} for $\Phi$.
The smallest $\lam$ satisfying eqn.(\ref{Lipschitz}) is called the
{\dfn $d$-Lipschitz constant} of $\Phi$.   More generally,
a {\dfn Lipschitz pseudometric} is a pseudometric
 $d:\sX^2\into\Real_+$ satisfying eqn.(\ref{Lipschitz}).

\example{\label{X:lipschitz0}
 Let $(\AV,\sX,\Phi)$ be a symbolic dynamical system, with network $(\dV,\connect)$.
For any $\fv\in\dV$ and $\lam>1$, we define
the pseudometric
$d_{\fv,\lam}:\AV\x\AV\into\Real_+$ as follows:  for all $\ba,\bb\in\AV$,
\beqn
\label{cantor.metric.on.digraph}
  d_{\fv,\lam}(\ba,\bb) \quad:=\quad \frac{1}{\lam^{R(\ba,\bb)}},
\quad\mbox{where}\quad
R(\ba,\bb)\quad:=\quad \max\set{r\in\Natur}{\ba_{\dB(\fv,r)} = \bb_{\dB(\fv,r)}}.
\eeqn
Then $d_{\fv,\lam}$ is a $\Phi$-Lipschitz pseudometric with constant $\lam$.
  To see this, let $\bx,\bx'\in\sX$.  If $R(\bx,\bx')=r$,
then $d(\bx,\bx')=\frac{1}{\lam^r}$.  But if $R(\bx,\bx')=r$, then
$R\lb[\Phi(\bx'),\Phi(\bx')\rb] \geq r-1$, so 
$d\lb[\Phi(\bx'),\Phi(\bx')\rb] \leq \frac{1}{\lam^{r-1}} = \lam\cdot d(\bx,\bx')$, as desired.
}

 The pseudometric $d_{\fv,\lam}$ in Example \ref{X:lipschitz0} is not
necessarily a true metric, unless $\Union_{r=0}^\oo \dB(\fv,r)$
$=\dV$, which is not the case unless $\fv$ is downstream from every
element of $\fv$.  For many digraphs, there is no vertex
with this property.  Instead, let us say that a subset $\dU\subset\dV$
is an {\dfn estuary} if, for every $\fv\in\dV$, there exists some $\fu\in\dU$
with $\fv\upstream \fu$.  For example, $\dV$ itself is an estuary.
If $(\dV,\connect)$ is biconnected,
then any nonempty subset of $\dV$ (even a singleton) is an estuary.
More generally, if $\dU$ contains at least one vertex from each biconnected
component of $\dV$, then $\dU$ is an estuary.

\example{\label{X:lipschitz}
Let $\dU\subset\dV$ be an estuary.  Let $\bc:=(c_\fu)_{\fu\in\dU}\in\Real_+^\dU$ be a sequence of
positive {\dfn coefficients} for the elements in $\dU$, such that
$\D \sum_{\fu\in\dU} c_\fu <\oo$ (this is possible because $\dU$ is
always countable, because $\dV$ is countable).
Fix $\lam>1$, and for all $\fu\in\dU$, let  $d_{\fu,\lam}$ be the
pseudometric from Example \ref{X:lipschitz0}.
Define the metric
$d_{\bc,\lam}:\sX\x\sX\into\Real_+$ by
\[
  d_{\bc,\lam}(\ba,\bb)\quad:=\quad \sum_{\fu\in\dU}^\oo c_\fu \, d_{\fu,\lam}(\ba,\bb).
\]
Then $d_{\bc,\lam}$ is a true metric, because $\Union_{r=0}^\oo \dB(\dU,r)=\dV$,
because $\dU$ is an estuary. Also, $d_{\bc,\lam}$ satisfies
eqn.(\ref{Lipschitz}), because
each of the pseudometrics $d_{\fu,\lam}$ satisfies eqn.(\ref{Lipschitz})
(by Example \ref{X:lipschitz0}).  Thus, $d_{\bc,\lam}$ is a Lipschitz metric for $\Phi$.
We say that  $d_{\bc,\lam}$ is {\dfn based at} $\dU$.
}

  Observe that the metric on $\AV$ in Example \ref{X:lipschitz} can be
defined for any digraph structure on $(\dV,\connect)$ (without
reference to any particular map $\Phi:\AV\into\AV$).  
Any Cantor dynamical system can be represented as a symbolic
dynamical system, so Example \ref{X:lipschitz} shows that any Cantor
dynamical system admits a Lipschitz metric.  Indeed, it admits many
such metrics, because the Lipschitz constant $\lam$, the estuary $\dU$,
and the coefficient system $\bc$ in Example \ref{X:lipschitz} can be chosen
arbitrarily.

\paragraph{Dimension and Entropy.}
Let $(\sX,d)$ be a metric space.  For any $\eps>0$,
a {\dfn open $\eps$-cover} is a covering of $\sX$ by open sets
whose diameters are each at most $\eps$.  Let $N_\eps(\sX)$ be the
minimal cardinality of any open $\eps$-cover of $(\sX,d)$.
We define 
\beqn
\begin{array}{rcl}
\label{metric.dimension.defn}
  \unddim(\sX,d)
&:=&
\D \liminf_{\eps\goto0} \ \frac{\log_\alp \lb[\log_\bet(N_\eps(\sX))\rb]}{\log_\alp\lb[-\log_\gam(\eps)\rb]}
\\ \\
\And \quad
  \bardim(\sX,d)
&:=&
\D \limsup_{\eps\goto0} \ \frac{\log_\alp \lb[\log_\bet(N_\eps(\sX))\rb]}{\log_\alp\lb[-\log_\gam(\eps)\rb]}.
\end{array}
\eeqn
(Here, $\alp,\bet,\gam>1$ are any constants, and need not be equal ---the limits in (\ref{metric.dimension.defn})
are independent of the choice of $\alp,\bet,\gam$).
If $\unddim(\sX,d)=\bardim(\sX,d)$, then we refer to their common value as
``$\dim(\sX,d)$'', the {\dfn dimension} of $(\sX,d)$.
Note that formula (\ref{metric.dimension.defn})
differs from the `box-counting dimension'
$\boxdim(\sX,d):=
\lim_{\eps\goto0} \ \frac{\log_2(N_\eps)}{-\log_2(\eps)}$ by the extra
logarithms.  Also,  $\dim(\sX,d)$ is meaningful
even when $\boxdim(\sX,d)=\oo$.
(Indeed, $\boxdim(\sX,d)$ is finite iff $\dim(\sX,d)=1$).  We will now show
that $\dim(\sX,d)$ is closely related to the `connectivity dimension'
of the digraph $(\dV,\connect)$.

Let $\unddim(\dV,\connect)$ and $\bardim(\dV,\connect)$ be as defined in
eqn.{\rm(\ref{connectivity.dimension})} of \S\ref{S:prelim}.   If $\dU\subseteq\dV$ is an estuary
for $\dV$, then Lemma \ref{constant.nhood.growth.rate}
implies that 
\beqn
\label{sup.limsup}
\bardim(\dV,\connect)\quad=\quad
\sup_{\fu\in\dU}\  \limsup_{r\goto\oo} \  \frac{\log\lb|\dB(\fu,r)\rb|}{\log(r)}.
\eeqn
We say $(\dV,\connect)$ has {\dfn uniform dimension} on $\dU$ if the `sup'
and `limsup' can be exchanged:
\beqn
\label{uniform.dimension}
   \bardim(\dV,\connect)\quad=\quad
  \limsup_{r\goto\oo} \ \sup_{\fu\in\dU}\  \frac{\log\lb|\dB(\fu,r)\rb|}{\log(r)}.
\eeqn
Heuristically, eqn.(\ref{uniform.dimension}) means that the limsups in eqn.(\ref{sup.limsup}) converge
`uniformly' on $\dU$.  For example, if $\dU$ is finite, then
$(\dV,\connect)$ automatically has uniform dimension on $\dU$. 

A nonnegative sequence  $\{c_j\}_{j=1}^\oo$ has {\dfn precipitous decay} if
\beqn
\label{precipitous}
\lim_{\eps\goto 0} \ \frac{\ln[J(\eps)]}{\ln\lb|\ln(\eps)\rb|}  \ \ = \ \  0,
\ \ \mbox{where, for all $\eps>0$,} \ \ 
J(\eps) \ \ := \ \  \min\set{J\in\Natur}{\sum_{j=J+1}^\oo c_j < \frac{\eps}{2}}.
\quad
\eeqn
Let $\dU\subseteq\dV$ and let $\bc=(c_\fu)_{\fu\in\dU}$ be some coefficient sequence.
Suppose we enumerate $\dU$ as $\dU=\{\fu_j\}_{j=0}^\oo$; then we can define
$c'_j:=c_{\fu_j}$ for all $j\in\Natur$;
then we say $\bc$ has {\dfn precipitous decay} if the sequence $\{c'_j\}_{j=1}^\oo$ 
has precipitous decay.

\example{(a) If $\{c_j\}_{j=1}^\oo$ has only finitely many nonzero terms,
then $\{c_j\}_{j=1}^\oo$ has  precipitous decay.
({\em Proof.} \  If $c_j=0$ for all $j\geq J_0$, then $J(\eps)\leq J_0$ for all $\eps$.)

(b) Let $c_j:=\exp(-e^j)\cdot e^j$ for all $j\in\Natur$.  Then 
 $\{c_j\}_{j=1}^\oo$ has  precipitous decay.
({\em Proof.} \ If $F(x):= -\exp(-e^x)$, then $F'(x) = \exp(-e^x)\cdot e^x$, so
$c_j = F'(j)$ for all $j\in\Natur$.  Thus,
$\sum_{j=J+1}^\oo c_j < \int_{J+1}^\oo F'(x) \ dx = -F(J+1)$, so
$J(\eps) \leq F^{-1}(-\eps/2) = \ln[-\ln(\eps/2)]$.)}

\Proposition{\label{entropy.vs.dimension}}
{
Let $(\dV,\connect)$ be a digraph,  let $\dU\subset\dV$ be an estuary, let $\lam>0$, and
let $d=d_{\lam,\bc}:\sX^2\into\Real_+$ be a
metric based on $\dU$, as in {\rm Example \ref{X:lipschitz}.}
Let $\sX\subseteq\AV$ be a pattern space.  
\bthmlist
\item Suppose $\dU':=\set{\fu\in\dU}{\undh_\fu(\sX)>0}$ is nonempty, and
let $\undD \D :=\sup_{\fu\in\dU'} \ \unddim_\fu(\dV,\connect)$.
Then $\unddim(\sX,d)\geq \undD$.
In particular, $\unddim(\sX,d)\geq \unddim(\dV,\connect)$.
\ethmlist
Let $\bc=(c_\fu)_{\fu\in\dU}$ be the coefficients
used to define $d$.
\bthmlist
\setcounter{enumii}{1}
\item  If  $\bc$ has {precipitous decay} and
 $(\dV,\connect)$ has {uniform dimension} on $\dU$,  then\linebreak 
${\bardim(\sX,d)\ \leq\  \bardim(\dV,\connect)}$.

\item In particular, if  $\dU$ is finite,  then 
$\bardim(\sX,d)\ \leq\  \bardim(\dV,\connect)$.
\ethmlist
}
\bthmprf 
Let $\dW\subset\dV$ be any finite set.
For all $\bw\in\sX_{\dW}$, let
$\inn{\bw}:=\set{\bx\in\sX}{\bx_{\dW}=\bw}$ be the
{\dfn cylinder set} defined by $\bw$.  The collection
$\sC_\dW:=\set{\inn{\bw}}{\bw\in\sX_\dW}$ is an open cover of $\sX$.

(a) \  Let $\unddel<\undD$.

\Claim{\label{entropy.vs.dimension.C2}
There exists $\eps_1>0$, $H>0$, and $L\in\Real$
 such that, for all $\eps\in\OO{0,\eps_1}$, we have
\[ \ln\lb(\log_2[N_\eps(\sX)]\rb) \quad>\quad
\ln(H) +  \unddel\cdot \ln(L-\log_\lam(\eps)).\]
\vspace{-2em}
}
\bclaimprf
For any $\eps>0$, let $\dU(\eps):=\set{\fu\in\dU}{c_\fu>\eps}$ 
(which is finite because $\bc$ is summable).
For all $\fu\in\dU(\eps)$, let $r_\fu(\eps):=\lfloor\log_\lam (c_\fu/\eps)\rfloor$.
Let ${\dW(\eps)}:=\D \Union_{\fu\in\dU(\eps)} \dB(\fu,r_\fu(\eps))$.

\subclaim{Let $\bx,\by\in\sX$. If $\bx_{\dW(\eps)}\neq \by_{\dW(\eps)}$, then $d(\bx,\by)>\eps$.\label{entropy.vs.dimension.C2.1}}
\bsubclaimprf
If $\bx_{\dW(\eps)}\neq \by_{\dW(\eps)}$, then there exists $\fu\in\dU(\eps)$ with
$\bx_{\dB(\fu,r_\fu(\eps))}\neq \by_{\dB(\fu,r_\fu(\eps))}$. Thus
\[
d(\bx,\by) \quad:=\quad \sum_{\fv\in\dU} c_\fv\, d_{\fv,\lam}(\bx,\by)
\quad \geq \quad 
c_\fu\cdot d_{\fu,\lam}(\bx,\by)
\quad \grt{(*)} \quad 
\frac{c_\fu}{\lam^{r_u(\eps)}}
\quad \geeeq{(\dagger)} \quad 
\frac{c_\fu\,\eps}{c_\fu} 
\quad=\quad \eps,
\]
as desired.  Here, $(*)$ is by eqn.(\ref{cantor.metric.on.digraph}),
and $(\dagger)$ is because $r_u(\eps)\leq \log_\lam (c_\fu/\eps)$.
\esubclaimprf

By hypothesis, there exists $\fu^*\in\dU'$ with
 $\unddim_{\fu^*}(\dB,\connect)> \unddel$.
Now, $\D\lim_{\eps\goto0}\ \dU(\eps)=\dU$ (because $c_\fu>0$ for all $\fu\in\dU$), so there exists $\eps_0>0$
such that, if $\eps\in\OO{0,\eps_0}$,
then  $\fu^*\in\dU(\eps)$.
Let $0<H<\undh_{\fu^*}(\sX)$.
Defining equations (\ref{connectivity.dimension}) and (\ref{entropy.defn}) in \S\ref{S:prelim}
say there exists $R>0$ such that, for all $r>R$, we have
$\D \frac{\ln\lb|\dB(\fu^*,r)\rb|}{\ln(r)}> \unddel$
and $\D\frac{\log_2\lb|\sX_{\dB(\fu^*,r)}\rb|}{\lb|\dB(\fu^*,r)\rb|} > H$.
But $\D\lim_{\eps\goto0}\ r_{\fu^*}(\eps)=\oo$.  Thus,
there exists $\eps_1\in\OO{0,\eps_0}$ such that, 
if $\eps\in\OO{0,\eps_1}$, then $r_{\fu^*}(\eps)>R$; hence
\begin{eqnarray}
\ln\lb|\dB(\fu,r_{\fu^*}(\eps))\rb|&> & \unddel\cdot \ln(r_{\fu^*}(\eps)) \label{entropy.vs.dimension.e1} \\
\And 
\log_2\lb|\sX_{\dB(\fu^*,r_{\fu^*}(\eps))}\rb| &>& H\cdot \lb|\dB(\fu^*,r_{\fu^*}(\eps))\rb|. \label{entropy.vs.dimension.e2}
\end{eqnarray}

Let $\eps\in\OO{0,\eps_1}$, and let $\sC_\eps$ be a minimal open $\eps$-cover;  then Claim \ref{entropy.vs.dimension.C2}.\ref{entropy.vs.dimension.C2.1} implies
that each cell of $\sC_\eps$ can intersect at most one cylinder
set from the cover $\sC_{\dW(\eps)}$.  Thus,
\beq
N_\eps(\sX) &=& |\sC_\eps| \quad\geq\quad |\sC_{\dW(\eps)}|
\quad=\quad |\sX_{\dW(\eps)}| 
\quad\geq\quad
\max_{\fu\in\dU(\eps)} \lb|\sX_{\dB(\fu,r_\fu(\eps))}\rb|
\\ &\geeeq{(*)}&
 \lb|\sX_{\dB(\fu^*,r_{\fu^*}(\eps))}\rb|.
\\
\mbox{Thus,}\quad
\log_2[N_\eps(\sX) ] &\geq &  
\log_2\lb|\sX_{\dB(\fu^*,r_{\fu^*}(\eps))}\rb|
\quad\grt{(\dagger)}\quad
 H \cdot\lb|\dB(\fu^*,r_{\fu^*}(\eps))\rb|\\
\mbox{Thus,}\quad
\ln\lb(\log_2[N_\eps(\sX)]\rb) &> &  
 \ln(H) + \ln\lb|\dB(\fu^*,r_{\fu^*}(\eps))\rb|
\\&\grt{(\ddagger)}&
  \ln(H) +  \unddel\cdot \ln(r_{\fu^*}(\eps))
\quad\grt{(\diamond)}\quad
  \ln(H) +  \unddel\cdot \ln(\log_\lam (c_{\fu^*}/\eps)-1)
\\&=&
  \ln(H) +  \unddel\cdot \ln(L-\log_\lam(\eps)),
\eeq
where $L:=\log(c_{\fu^*})-1$.
Here, $(*)$ is because $\fu^*\in\dU(\eps)$ because $\eps<\eps_0$.
$(\dagger)$ is by (\ref{entropy.vs.dimension.e2}), \ 
$(\ddagger)$ is by (\ref{entropy.vs.dimension.e1}), \ and $(\diamond)$ is because 
 $r_{\fu^*}(\eps):=\lfloor\log_\lam (c_{\fu^*}/\eps)\rfloor>\log_\lam (c_{\fu^*}/\eps)-1$.
\eclaimprf
\beq
\mbox{We now have}\qquad
 \unddim(\sX,d)
&:\eeequals{(*)}&
\liminf_{\eps\goto0} \ \frac{\ln \lb[\log_2(N_\eps(\sX))\rb]}{\ln\lb[-\log_\lam(\eps)\rb]}
\\& \geeeq{(\dagger)} &
\liminf_{\eps\goto0} \  \frac{ \ln(H) +  \unddel\cdot \ln(L-\log_\lam(\eps)) }{\ln\lb[-\log_\lam(\eps)\rb]}
\quad=\quad \unddel,
\eeq
where $(*)$ is by setting $\alp:=e$, \ $\bet:=2$, and $\gam:=\lam$ in definition (\ref{metric.dimension.defn}), while
$(\dagger)$ is by Claim \ref{entropy.vs.dimension.C2}.  This holds for any $\unddel<\undD$.  Thus, we conclude that
$\unddim(\sX,d)\geq \undD$, as desired.

(b) \  Fix some enumeration $\dU=\{\fu_j\}_{j=0}^\oo$ and
define $c_j:=c_{\fu_j}$ for all $j\in\Natur$.
For all $\eps>0$, let $\D J(\eps)$ be as in eqn.(\ref{precipitous}).
Let $\bardel>\bardim(\dV,\connect)$.  Let $\D S:=\sum_{j=0}^\oo c_j<\oo$. 

\Claim{\label{entropy.vs.dimension.C3}
There exists $\eps_1>0$ and constants $L_1,L_2>0$ such that, for any $\eps\in\OO{0,\eps_1}$:
\[
\ln\lb(\log_2[N_\eps(\sX)]\rb)
\quad\leq\quad 
\ln[J(\eps)]\ + \  L_1 \ + \   \bardel\cdot \ln[L_2 -\log_\lam(\eps)].
\]
\vspace{-2em}
}
\bclaimprf
For any $\eps>0$, let $r(\eps):=\lceil \log_\lam(2S/\eps)\rceil$,
and let $\D {\dW(\eps)}:=\Union_{j=0}^{J(\eps)} \dB[\fv_j,r(\eps)]$.

\subclaim{Let $\bx,\by\in\sX$.
If $\bx_{\dW(\eps)}=\by_{\dW(\eps)}$, then $d(\bx,\by)<\eps$.\label{entropy.vs.dimension.C3.1}}
\bsubclaimprf
We have
\beq
d(\bx,\by) &:=& \sum_{j=0}^\oo c_j\, d_{\fv_j,\lam}(\bx,\by)
\quad=\quad
\sum_{j=0}^{J(\eps)} c_j\, d_{\fv_j,\lam}(\bx,\by) \ + \ \sum_{j=J(\eps)+1}^\oo c_j\, d_{\fv_j,\lam}(\bx,\by)
\\ &\leq&
\lb(\sum_{j=0}^{J(\eps)} c_j\rb) \cdot \max_{0\leq j\leq J(\eps)} \lb(d_{\fv_j,\lam}(\bx,\by)\rb)
 \ + \  \lb(\sum_{j=J(\eps)+1}^\oo c_j\rb) \cdot \max_{j\geq J(\eps)} \lb(d_{\fv_j,\lam}(\bx,\by)\rb)
\\&\leeeq{(*)} &
\frac{S}{\lam^{r(\eps)}} \ + \ \frac{\eps}{2} \cdot 1
\quad\leeeq{(\dagger)}\quad
\frac{\eps S}{2S} \ + \ \frac{\eps}{2} 
\quad=\quad
\frac{\eps}{2}+\frac{\eps}{2} \quad=\quad \eps, \qquad\mbox{as desired.}
\eeq
 Here, $(*)$ is because $d_{\fv_j,\lam}(\bx,\by)\leq\frac{1}{\lam^{r(\eps)}}$ for all $j\in\CC{0...J(\eps)}$
because 
$\bx_{\dB[\fv_j,r(\eps)]}=\bx_{\dB[\fv_j,r(\eps)]}$ for all $j\in\CC{0...J(\eps)}$;
meanwhile, $\D\sum_{j=0}^{J(\eps)} c_j\leq \sum_{j=0}^\oo c_j = S$, and
$\D\sum_{j=J(\eps)+1}^\oo c_j<\frac{\eps}{2}$ by definition of $J(\eps)$.
Finally,  $(\dagger)$ is because $r=\lceil \log_\lam(2S/\eps)\rceil$.
\esubclaimprf

Now, $\bardel>\bardim(\dV,\connect)$, so eqn.(\ref{uniform.dimension}) yields
some $R\in\Natur$ such that, for all $r\in\Natur$:
\beqn
\label{entropy.vs.dimension.e3}
\mbox{if $r>R$, then}\qquad
 \sup_{\fu\in\dU} \ \ln\lb|\dB(\fu,r)\rb| \quad<\quad \bardel\cdot \ln(r).
\eeqn
Now, $\D \lim_{\eps\goto0}\ r(\eps)=\oo$, so there exists 
$\eps_1>0$ such that, if $\eps\in\OO{0,\eps_1}$, then  $r(\eps)>R$.
Let $\eps\in\OO{0,\eps_1}$.  Claim \ref{entropy.vs.dimension.C3}.\ref{entropy.vs.dimension.C3.1} implies that $\sC_{\dW(\eps)}$ is an $\eps$-open cover of $\sX$.  Thus,
\beq
\nonumber
N_\eps(\sX)&\leq& |\sC_{\dW(\eps)}|\quad=\quad |\sX_{\dW(\eps)}|\quad\leq\quad \prod_{j=0}^{J(\eps)} \lb|\sX_{\dB(\fv_j,r(\eps))}\rb|.\\
\nonumber
\mbox{Thus,}\quad
\log_2[N_\eps(\sX)] &\leq&
\sum_{j=0}^{J(\eps)} \log_2\lb|\sX_{\dB(\fv_j,r(\eps))}\rb| 
\quad\leq\quad
\sum_{j=0}^{J(\eps)} \lb(\log_2|\sA|\rb)\cdot\ \lb|\dB(\fv_j,r(\eps))\rb| 
\\&\leq &
J(\eps)\cdot \log_2|\sA|\cdot \max_{0\leq j\leq J(\eps)}  \lb|\dB(\fv_j,r(\eps))\rb|
\nonumber \\&\leq &
J(\eps)\cdot \log_2|\sA|\cdot \sup_{\fu\in\dU} \ \lb|\dB(\fu,r(\eps))\rb|.\\
\mbox{Thus,}\quad
\ln\lb(\log_2[N_\eps(\sX)]\rb) 
&\leq&
\ln[J(\eps)]\ + \ \ln(\log_2|\sA|) \ + \   \sup_{\fu\in\dU} \ \ln\lb|\dB(\fu,r(\eps))\rb|\\
&\leeeq{(*)}&
\ln[J(\eps)]\ + \  \ln(\log_2|\sA|)  \ + \   \bardel\cdot \ln(r(\eps)) \\
&\leeeq{(\dagger)}&
\ln[J(\eps)]\ + \  L_1 \ + \   \bardel\cdot \ln(r(\eps)) \\
&\leeeq{(\ddagger)}&
\ln[J(\eps)]\ + \  L_1 \ + \   \bardel\cdot \ln[1 + \log_\lam(2S/\eps)] \\
&\leq&
\ln[J(\eps)]\ + \  L_1 \ + \   \bardel\cdot \ln[L_2 -\log_\lam(\eps)],
\eeq
where $L_2:=1+\log_\lam(2S)$.  Here,
$(*)$ is by eqn.(\ref{entropy.vs.dimension.e3}), because $r(\eps)>R$ because $\eps\in\OO{0,\eps_1}$.  In
$(\dagger)$ we define  $L_1:=\ln(\log_2|\sA|)$].  Finally,
$(\ddagger)$ is because  $r:=\lceil \log_\lam(2S/\eps)\rceil
\leq 1 +\log_\lam(2S/\eps)$.
\eclaimprf
\beq
\mbox{We now have}\qquad
  \bardim(\sX,d)
&:\eeequals{(*)}&
\limsup_{\eps\goto0} \ \frac{\ln \lb[\log_2(N_\eps(\sX))\rb]}{\ln\lb[-\log_\lam(\eps)\rb]}
\\
&\leeeq{(\dagger)}&
\limsup_{\eps\goto0} \ \frac{\ln[J(\eps)]\ + \  L_1 \ + \   \bardel\cdot \ln[L_2 -\log_\lam(\eps)]}{\ln\lb[-\log_\lam(\eps)\rb]} 
\\&=&
\bardel \ + \ \lim_{\eps\goto0} \ \frac{\ln[J(\eps)]}{\ln\lb[-\log_\lam(\eps)\rb]} 
\quad\eeequals{(\ddagger)}\quad \bardel. 
\eeq
Here, $(*)$ is by setting $\alp:=e$, \ $\bet:=2$, and $\gam:=\lam$ in  definition (\ref{metric.dimension.defn}), and
$(\dagger)$ is by Claim 2.  Meanwhile, $(\ddagger)$ is because $\{c_j\}_{j=1}^\oo$ has
precipitous decay.  

Thus works for any $\bardel>\bardim(\dX,\connect)$;
we conclude that $ \bardim(\sX,d)\leq  \bardim(\dX,\connect)$.

(c) follows immediately from (b), because if $\dU$ is finite, then clearly
$(\dV,\connect)$ has uniform dimension on $\dU$, and $\bc$  has
precipitous decay.  
\ethmprf

Let $(\dV,\connect)$ be a dimensionally homogeneous digraph
[i.e. $\unddim(\dV,\connect)=\bardim(\dV,\connect)$].
If $\sX\subseteq\AV$ is a pattern space, and
$d:\sX^2\into\Real_+$ is a Cantor metric, then we say that $d$ is 
{\dfn dimensionally compatible} if $\dim(\sX,d)=\dim(\dV,\connect)$.
Proposition \ref{entropy.vs.dimension} suggests that for `most' 
dimensionally homogeneous digraphs, any pattern space with nonzero
entropy admits a dimensionally compatible
metric.  In light of Example \ref{X:lipschitz}, this means that
`most' symbolic dynamical systems admit dimensionally compatible
Lipschitz metrics.  For example, we have the following result:

\Corollary{\label{dimensionally.compatible}}
{
Let  $(\dV,\connect)$ be a dimensionally homogeneous digraph with
a finite estuary $\dU$ {\rm (e.g. a biconnected digraph)}.  There exists a metric $d$ on $\AV$ such that, if
$\sX\subseteq\AV$ is any pattern space with $\undh_\fu(\sX)>0$ for some $\fu\in\dU$,
then $\dim(\sX,d)=\dim(\dV,\connect)$.  Furthermore, if $\Phi:\AV\into\AV$
is a continuous map with network $(\connect)$, then $\Phi$
is $d$-Lipschitz.
}
\bthmprf Let  $d$ be a metric based
on $\dU$, as in Example \ref{X:lipschitz}.  Then
$\bardim(\sX,d)\ \leeeq{(*)}\  {\bardim(\dV,\connect)}$ \ $ \eeequals{(\dagger)}\
 {\dim(\dV,\connect)}$,
where $(*)$ is by Proposition \ref{entropy.vs.dimension}(c) and $(\dagger)$ is by dimensional
homogeneity.
On the other hand, 
$\unddim(\sX,d) \ \geeeq{(*)} \  \unddim(\dV,\connect) \ \eeequals{(\dagger)}\  \dim(\dV,\connect)$,
where $(*)$ is by Proposition \ref{entropy.vs.dimension}(a) and
$(\dagger)$ is by dimensional
homogeneity.  We conclude that $\dim(\dV,\connect)\leq \unddim(\sX,d)\leq \bardim(\sX,d)\leq \dim(\dV,\connect)$;
hence  $\dim(\sX,d)$ is well-defined and $\dim(\sX,d)=\dim(\dV,\connect)$.
The fact that $d$ is $\Phi$-lipschitz was demonstrated in  Example \ref{X:lipschitz}.
\ethmprf

\example{\label{X:cantor.dimension}  
(a)   Let $\dV=\Zahl^2$ have the Cayley digraph structure
induced by generating set $\dB:=\{(\pm1,0)\}$, $\{(0,\pm1)\}$.
Then $(\Zahl^2,\connect)$ is biconnected, so any singleton set is an estuary.
So, let $\fo=(0,0)$ be the origin, 
and let $\dU:=\{\fo\}$;  set $c_\fo=1$ and $c_\fz=0$ for all nonzero $\fz\in\Zahl^2$.

For any $r>0$, we have $\dB(\fo,r):=\set{\fz\in\Zahl^2}{ |z_1|+|z_2|\leq r}$.
Let $\lam=2$; then the  metric $d_{\bc,\lam}$ from Example \ref{X:lipschitz}
becomes the standard Cantor metric on $\sA^{\Zahl^2}$:
\[
  d(\ba,\ba') \quad:=\quad \frac{1}{2^R},\quad\mbox{where}\quad
R \ := \ \max\set{r\in\Natur}{\ba_{\dB(\fo,r)} = \ba'_{\dB(\fo,r)}}.
\]
If $\sX\subseteq\sA^{\Zahl^2}$ is any subshift with positive topological entropy,
then Proposition \ref{entropy.vs.dimension} says
 ${\dim(\sX,d)}={\dim(\Zahl^2,\connect)}=2$.  If $\Phi:\sA^{\Zahl^2}\into\sA^{\Zahl^2}$
is any cellular automaton whose local rule has neighbourhood $\{0\}\disj\dB$,
then $\Phi$ is $d$-Lipschitz, by Example \ref{X:lipschitz}.

(b) \ By a similar argument, if $\dG$ is any finitely generated group with growth dimension
$D$ and a biconnected Cayley digraph structure, then we can
construct a Cantor metric $d$ on $\sA^\dG$ such that,
if $\sX\subseteq\sA^\dG$ is any subshift with positive topological entropy, then $\dim(\sX,d)=D$.
Furthermore, if $\Phi$ is any CA on $\sX$, we can design $d$
to be $\Phi$-Lipschitz.

(c) \ However, it is possible to construct zero-entropy subshifts of $\AZD$
with dimensions less than $D$.  For example, treat $\sA^{\Zahl^2} \cong
\prod_{z\in\Zahl} \AZ$ in the obvious way, so that any $\ba\in\sA^{\Zahl^2}$
has the form $\ba=(\ldots,\ba_{-1},\ba_0,\ba_1,\ba_2,\ldots)$, where
$\ba_z\in\AZ$ for all $z\in\Zahl$.  Let
$\sX:=\set{ (\ldots,\ba,\ba,\ba,\ldots)}{\ba\in\AZ}\subset\sA^{\Zahl^2}$.
Then $\htop(\sX)=0$ and  $\dim(\sX,d)=1$.}

Let $(\sX,d)$ and $(\sX',d')$ be two metric spaces.
A continuous function $\Gam:\sX\into\sX'$ is {\dfn $(d,d')$-\Holder}
if, there exist $\eta,\lam\in\OO{0,\oo}$ such that,
for any $\bx_1,\bx_2\in\sX$,
\beqn
\label{Holder}
  d'\lb(\Phi(\bx_1), \ \Phi(\bx_2)\rb)
\quad\leq\quad
 \lam\cdot d(\bx_1,\bx_2)^\eta.
\eeqn
For example, any Lipschitz function is \Holder, with $\eta= 1$. 
If $\Gam$ is a homeomorphism, then we say $\Gam$ is {\dfn $(d,d')$-bi\Holder} if
both $\Gam$ and $\Gam^{-1}$ are \Holder\  (possibly with
different values of $\eta$ and/or $\lam$).

\Proposition{\label{dimension.holder.invariant}}
{
Let $(\sX,d)$ and $(\sX',d')$ be  metric spaces.
\bthmlist
\item Let $\Gam:\sX\into\sX'$ be a $(d,d')$-\Holder\  surjection.  Then 
$\dim(\sX,d)\geq \dim(\sX',d')$.

\item If $\Gam$ is a $(d,d')$-bi\Holder\  homeomorphism, then
$\dim(\sX,d)= \dim(\sX',d')$.
\ethmlist
}
\bthmprf (b) follows from (a).  To see
(a), suppose $\Gam$ is $(d,d'$)-\Holder, and let $\eta,\lam>0$ be as in
eqn.(\ref{Holder}).

\Claim{For any $\eps>0$, \ $N_{\eps}(\sX)  \ \geq \  N_{\lam\eps^\eta}(\sX')$.}
\bclaimprf 
Let $\sO:=\{\bO_1,\bO_2,\ldots,\bO_N\}$ be any open $\eps$-cover
of $\sX$.  Then for each $n\in\CC{1...N}$, \ the set $\Phi(\bO_n)$ is
open (because $\Phi$ is an open map, being
a continuous
surjection onto a compact space), and has diameter at most $\lam\eps^\eta$
by eqn.(\ref{Holder}).   The collection
$\Phi(\sO):=\{\Phi(\bO_1),\ldots,\Phi(\bO_n)\}$ together covers $\sX'$,
because $\sO$ covers $\sX$ and  $\Phi$ is surjective. Thus, $\Phi(\sO)$
is a $(\lam\eps^\eta)$-diameter open
cover of $\sX'$.  
  If $\sO$ is a minimal open $\eps$-cover of
$\sX$, then $N_{\eps}(\sX)=N$.  Since $\Phi(\sO)$ is a $(\lam\eps^\eta)$-cover
of $\sX'$ with $N$ pieces, we have $N_{\lam\eps^\eta}(\sX')\leq N$.
\eclaimprf

It follows that
\beq
\dim(\sX,d)
&=&
\lim_{\eps\goto0} \ \frac{\log[\log(N_\eps(\sX))]}{\log[-\log(\eps)]}
\quad\geeeq{(\ddagger)}\quad
\lim_{\eps\goto0} \ \frac{\log[\log(N_{\lam\eps^\eta}(\sX'))]}{\log[-\log(\eps)]}
\\&=&
\lim_{\eps\goto0} \ \lb(\frac{\log[-\log(\lam\eps^\eta)]}{\log[-\log(\eps)]}\rb)
 \cdot
\lb( \frac{\log[\log(N_{\lam\eps^\eta}(\sX'))]}{\log[-\log(\lam\eps^\eta)]}\rb)
\\
&\eeequals{(*)}&
\lim_{\eps\goto0}
\ \frac{\log[\log(N_{\lam\eps^\eta}(\sX'))]}{\log[-\log(\lam\eps^\eta)]}
\quad\eeequals{(\dagger)}\quad
\lim_{\eps'\goto0} \ 
\frac{\log[\log(N_{\eps'}(\sX'))]}{\log[-\log(\eps')]}
\quad=\quad
\dim(\sX',d).
\eeq
Here, $(\ddagger)$ is by Claim 1, and 
$(\dagger)$ is
where we make the change of variables $\eps'=\lam\eps^\eta$.
Finally, $(*)$ is because
\beq
\lim_{\eps\goto0} \ \frac{\log[-\log(\lam\eps^\eta)]}{\log[-\log(\eps)]}
&=&
 \lim_{\eps\goto0} \ \frac{\log[-\log(\lam)-\eta\log(\eps)]}{\log[-\log(\eps)]}
\quad\eeequals{(H)}\quad
 \lim_{\eps\goto0} \ \frac{\frac{-\eta/\eps}{-\log(\lam)-\eta\log(\eps)}}
{\frac{-1/\eps}{-\log(\eps)}}
\\&=&
 \lim_{\eps\goto0} \frac{\eta\log(\eps)}{\log(\lam)+\eta\log(\eps)}
\quad=\quad
1,
\eeq
where (H) is by L'Hospital's rule.
\ethmprf

It follows that the connectivity network dimension of a symbolic dynamical
system is invariant under bi\Holder\  topological conjugacy.

\Corollary{\label{network.dimension.invariant.under.holder.conjugacy}}
{
Let $(\AV,\sX_1,\Phi_1)$ and $(\BW,\sX_2,\Phi_2)$ be two symbolic dynamical systems,
and let $d_1$ and $d_2$ be
dimensionally compatible Lipschitz metrics 
on $\sX_1$ and $\sX_2$ respectively {\rm(e.g. as given by Corollary \ref{dimensionally.compatible})}.
\bthmlist
\item If there is a factor mapping $(\sX_1,\Phi_1)\into(\sX_2,\Phi_2)$
which is $(d_1,d_2)$-\Holder, 
then\linebreak {$\dim(\dV,\connect_{\!\!_1})\geq \dim(\dW,\connect_{\!\!_2})$}.

\item If $(\sX_1,\Phi_1)$ and $(\sX_2,\Phi_2)$ are conjugate via a bi-\Holder\  homeomorphism,
then \linebreak $\dim(\dV,\connect_{\!\!_1})= \dim(\dW,\connect_{\!\!_2})$.
\ethmlist
}

\Remark{Clearly, a continuous function $\Phi:\AV\into\AV$ also admits other
Lipschitz metrics which are {\em not} dimensionally compatible.
For example, let $\dU=\dV$ in Example \ref{X:lipschitz}.
If the coefficient system $\bc$ decays slowly enough, we can
make $\dim(\AV,d_{\bc,\lam})$ arbitrarily large. However, if $\undh(\sX)>0$, then
Proposition \ref{entropy.vs.dimension}(a) says 
it is not possible to make $\dim(\AV,d_{\bc,\lam})$
smaller than $\unddim(\dV,\connect)$ for any choice of $\bc$.
  }

\section*{Conclusion}

For any symbolic dynamical system $(\AV,\sX,\Phi)$, one can define a
digraph structure $(\connect)$ on $\dV$.   We have shown that certain
topological-dynamical properties of $(\AV,\sX,\Phi)$ are related to
the connectivity $(\dV,\connect)$, and in particular, to its dimension.
What other dynamical properties of $(\AV,\sX,\Phi)$ are influenced by
the geometry of $(\dV,\connect)$?

One could also go the other way.  Starting with an infinite digraph
$(\dV,\connect)$, consider a randomly generated self-map
$\Phi:\AV\into\AV$, such that $(\connect)$ is the network of $\Phi$.
What are the `generic' (i.e. almost-certain) properties of
$(\AV,\Phi)$, and how do they depend on the geometry of
$(\dV,\connect)$?
  
  For example, \S\ref{S:propagation} suggests the following
conjecture: {\em If $\dim(\dV,\connect)\leq 1$, then almost surely,
$(\AV,\Phi)$ is equicontinuous.  If $\dim(\dV,\connect)> 1$, then
almost surely, $(\AV,\Phi)$ is sensitive.}  (The intuition here comes
from percolation theory).  However, Figure \ref{fig:expansive.2d}
shows that something more than dimension is required; this network has
dimension 2, but it has an infinite number of cut points, so a random
mapping $\Phi$ with this network is almost-surely equicontinuous.
Thus, the conjecture above must be augmented with some kind of `regularity'
condition on $(\dV,\connect)$.

  A closely related question: Suppose we take a system
$(\AV,\Phi)$ and `mutate' it, by changing the local rule at a small
number of vertices.  What topological-dynamical properties are
`robust' under such mutations, and how does this depend on the
geometry of $(\dV,\connect)$?

\paragraph{Acknowledgments.}  
This research began during a research leave at Wesleyan University,
and was partially supported by the Van Vleck Fund;  I am grateful to
Ethan Coven for his generous hospitality.  This work benefitted from
conversations with Ethan, as well as Pierre Tisseur and Reem Yassawi.
This research was also supported by NSERC Grant \#262620-2008.

{\footnotesize
\bibliographystyle{alpha}
\bibliography{bibliography}
}

\end{document}